\documentclass[11pt, twoside, leqno, amsmath, amsfonts,
eucal, mathscr]{aomamlt2e} \pageno{1} \received{June 30, 2006}
\revised{}

\newtheorem{theorem}{Theorem}[section]

\newtheorem{lemma}[theorem]{Lemma}
\newtheorem{proposition}[theorem]{Proposition}
\newtheorem{corollary}[theorem]{Corollary}
\newtheorem{conjecture}[theorem]{Conjecture}

\theoremstyle{plain}
\newtheorem{remark}[theorem]{Remark}
\newtheorem{example}[theorem]{Example}

\numberwithin{equation}{section}

\renewcommand{\theequation}{\thesubsection.\arabic{equation}}

\newcommand{\M}{\mathcal{M}}
\newcommand{\T}{\mathscr{F}}
\newcommand{\Mbar}{\overline{M}}
\newcommand{\ft}{ft}
\newcommand{\Pd}{P_d}
\newcommand{\Pbeta}{P_\beta}
\newcommand{\q}{q}
\newcommand{\n}{\pi}

\newcommand{\e}{e}
\newcommand{\s}{\varepsilon}

\begin{document}
\currannalsline{0}{2007}

\title{Flops, motives and\\ invariance of quantum rings}

\acknowledgements{}
\author{Yuan-Pin Lee, Hui-Wen Lin and Chin-Lung Wang}
\institution{Department of Mathematics, University
of Utah, Salt Lake City, Utah 84103, U.S.A.\\
\email{yplee@math.utah.edu}\\ \\
Department of Mathematics, National Central University, Chung-Li,
Taiwan.\\ \email{linhw@math.ncu.edu.tw}\\ \\ Department of
Mathematics, National Central University, Chung-Li, Taiwan.
National Center for Theoretic Sciences, Hsinchu, Taiwan. \\
\email{dragon@math.ncu.edu.tw; dragon@math.cts.nthu.edu.tw}}

\shorttitle{}


\begin{abstract}
For ordinary flops, the correspondence defined by the graph
closure is shown to give equivalence of Chow motives and to
preserve the Poincar\'e pairing. In the case of simple ordinary
flops, this correspondence preserves the big quantum cohomology
ring after an analytic continuation over the extended K\"ahler
moduli space.

For Mukai flops, it is shown that the birational map for the local
models is deformation equivalent to isomorphisms. This implies
that the birational map induces isomorphisms on the full quantum
rings and all the quantum corrections attached to the extremal ray
vanish.
\end{abstract}

\setcounter{section}{-1}

\section{Introduction} \label{s:0}

\subsection{Statement of main results} \label{s:0.1}

Let $X$ be a smooth complex projective manifold and $\psi:X\to
\bar X$ a flopping contraction in the sense of minimal model
theory, with $\bar\psi: Z\to S$ the restriction map on the
exceptional loci. Assume that
\begin{itemize}
\item[(i)] $\bar\psi$ equips $Z$ with a $\mathbb{P}^r$-bundle
structure $\bar\psi : Z = \mathbb{P}_S(F) \rightarrow S$ for some
rank $r + 1$ vector bundle $F$ over a smooth base $S$, \item[(ii)]
$N_{Z/X}|_{Z_s} \cong \mathscr{O}_{\mathbb{P}^r}(-1)^{\oplus(r +
1)}$ for each $\bar\psi$-fiber $Z_s$, $s\in S$.
\end{itemize}
It is not hard to see that the corresponding {\it ordinary
$\mathbb{P}^r$ flop} $f: X \dashrightarrow X'$ exists. An ordinary
flop is called \emph{simple} if $S$ is a point.

For a $\mathbb{P}^r$ flop $f: X \dashrightarrow X'$, the graph
closure $[\bar\Gamma_f] \in A^*(X \times X')$ identifies the Chow
motives $\hat{X}$ of $X$ and $\hat{X'}$ of $X'$. Indeed, let $\T
:= [\bar\Gamma_f]$ then the transpose $\T^*$ is
$[\bar\Gamma_{f^{-1}}]$. One has the following theorem.

\begin{theorem} \label{t:1.1}
For an ordinary $\mathbb{P}^r$ flop $f:X\dashrightarrow X'$, the
graph closure $\T := [\bar\Gamma_f]$ induces $\hat X \cong \hat
X'$ via $\T^*\circ \T = \Delta_X$ and $\T\circ \T^* =
\Delta_{X'}$. In particular, $\T$ preserves the Poincar\'e pairing
on cohomology groups.
\end{theorem}

While the ring structure is in general not preserved under $\T$,
the quantum cohomology ring is, when the analytic continuation on
the Novikov variables is allowed.

\begin{theorem} \label{t:1.2} 
The big quantum cohomology ring is invariant under simple ordinary
flops, after an analytic continuation over the extended K\"ahler
moduli space.
\end{theorem}

A contraction $(\psi, \bar\psi): (X,Z) \to (\bar X, S)$ is of
\emph{Mukai type} if $Z = \mathbb{P}_S(F) \rightarrow S$ is a
projective bundle under $\bar\psi$ and $N_{Z/X} = T^*_{Z/S}$. The
corresponding algebraic flop $f: X \dashrightarrow X'$ exists and
its local model can be realized as a {\it slice} of an ordinary
flop. The following result is proved based upon our understanding
of local geometry of Mukai flops.

\begin{theorem} \label{t:1.3}
Let $f: X \dashrightarrow X'$ be a Mukai flop. Then $X$ and $X'$
are diffeomorphic, and have isomorphic Hodge structures and full
Gromov--Witten theory. In fact, any local Mukai flop is a limit of
isomorphisms and all quantum corrections attached to the extremal
ray vanish.
\end{theorem}

\subsection{Motivations} \label{1.2}

This paper is the first of our study of the relationship between
birational geometry and Gromov--Witten theory. Our motivations
come from both fields.

\subsubsection*{$K$-equivalence in birational geometry}
Two ($\mathbb{Q}$-Gorenstein) varieties $X$ and $X'$ are
\emph{$K$-equivalent} if there exist birational morphisms $\phi: Y
\to X$ and $\phi': Y \to X'$ with $Y$ smooth such that
\[
  \phi^* K_X = {\phi'}^* K_{X'}.
\]
$K$-equivalent smooth varieties have the same Betti numbers
(\cite{vB} \cite{Wang1}, see also \cite{Wang2} for a survey on
recent development). However, the cohomology ring structures are
in general different. Two natural questions arise here:

\begin{enumerate}
\item Is there a \emph{canonical} correspondence between the
cohomology groups of $K$-equivalent smooth varieties?

\item Is there a modified ring structure which is invariant under
the $K$-equivalence relation?
\end{enumerate}

The following conjecture was advanced by Y.~Ruan \cite{yR2} and
the third author \cite{Wang2} in response to these questions.
\begin{conjecture} \label{c:1.4}
$K$-equivalent smooth varieties have canonically isomorphic
quantum cohomology rings over the extended K\"ahler moduli spaces.
\end{conjecture}

The choice to start with ordinary flops is almost obvious.
Ordinary flops are not only the first examples of $K$-equivalent
maps, but also crucial to the general theory. In fact, one of the
goals of this paper is to study some of their fundamental
properties.

\subsubsection*{Functoriality in Gromov--Witten theory}
In the Gromov--Witten theory, one is led to consider the problem
of functoriality in quantum cohomology. Quantum cohomology is not
functorial with respect to the usual operations: pull-backs,
push-forwards, etc.. Y.~Ruan \cite{yR1} has proposed to study the
\emph{Quantum Naturality Problem}: finding the ``morphisms'' in
the ``category'' of symplectic manifolds for which the quantum
cohomology is ``natural''.

The main reason for lack of functoriality comes from the dimension
count of the moduli of stable maps, where Gromov--Witten
invariants are defined. (See \S3.1 for the relevant definitions.)
For example, given a birational morphism $f: Y \to X$, there is an
induced morphism from moduli of maps to $Y$ to moduli of maps to
$X$. However, the (virtual) dimensions of the two moduli spaces
are equal only if $Y$ and $X$ are $K$-equivalent. When the virtual
dimensions of moduli spaces are different, the non-zero integral
on moduli space of maps to $X$ will be ``pulled-back'' to a zero
integral on moduli space of maps to $Y$. Therefore,
$K$-equivalence appears to be a necessary condition for this type
of functoriality. Conjecture~\ref{c:1.4} suggests that the
$K$-equivalence is also sufficient. We note here that there is of
course no $K$-equivalent morphism between smooth varieties and a
``flop-type'' transformation is needed.

Theorem~\ref{t:1.2} can therefore be considered as establishing
some functoriality of the genus zero Gromov--Witten theory in this
direction. The higher genus case will be discussed in a separate
paper.

\subsubsection*{Crepant resolution conjecture}
Conjecture~\ref{c:1.4} can also be interpreted as a consistency
check for the \emph{Crepant Resolution Conjecture} \cite{yR2}
\cite{BG}. In general, there are more than one possible crepant
resolution, but different crepant resolutions are $K$-equivalent.
The consistency check naturally leads to a special version of
Conjecture~\ref{c:1.4}.

\subsection{Contents of the paper} \label{1.3}
\S\ref{s:1} studies the geometry of ordinary flops. The existence
of ordinary flops is proved and explicit description of local
models is given.

\S\ref{s:2} is devoted to the correspondences and Chow motives of
projective smooth varieties under an ordinary flop. The main
result of this section is Theorem~\ref{t:1.1} alluded above. The
ring structure is, however, not preserved. For a simple
$\mathbb{P}^r$-flop, let $h$ be the hyperplane class of $Z =
\mathbb{P}^r$ and let $\alpha_{i} \in H^{2 l_i}(X)$, with $l_i \le
r$ and $l_1 + l_2 + l_3 = \dim X = 2r + 1$.
\begin{proposition}\label{classical-0.5}
\[
(\T\alpha_1.\T\alpha_2.\T\alpha_3) = (\alpha_1.\alpha_2.\alpha_3)
+ (-1)^r (\alpha_1.h^{r - l_1})(\alpha_2.h^{r -
l_2})(\alpha_3.h^{r - l_3}).
\]
\end{proposition}

For Calabi-Yau threefolds under a simple $\mathbb{P}^1$ flop, it
is well known in the context of string theory (see e.g.\
\cite{Witten}) that the defect of the classical product is exactly
remedied by the quantum corrections attached to the extremal rays.
This picture also emerged as part of Morrison's cone conjecture on
birational Calabi-Yau threefolds \cite{Morrison} where
Conjecture~\ref{c:1.4} for Calabi-Yau threefolds was proposed. For
threefolds Conjecture~\ref{c:1.4} was proved by A.~Li and Y.~Ruan
\cite{LiRu}. Their proof has three ingredients:
\begin{enumerate}
\item A symplectic deformation and decomposition of $K$-equivalent
maps into composite of ordinary $\mathbb{P}^1$ flops,

\item the multiple cover formula for $\mathbb{P}^1\cong C\subset
X$ with $N_{C/X}\cong \mathscr{O}(-1)^{\oplus 2}$, and their main
contribution:

\item the theory of relative Gromov-Witten invariants and the
degeneration formula.
\end{enumerate}

In \S\ref{s:3} a higher dimensional version of ingredient (2) is
proved:

\begin{theorem} \label{t:1.7} 
Let $Z = \mathbb{P}^r \subset X$ with $N_{Z/X} \cong
\mathscr{O}(-1)^{r + 1}$. Let $\ell$ be the line class in $Z$.
Then for all $\alpha_i\in H^{2l_i}(X)$ with $1\le l_i \le r$,
$\sum_{i = 1}^n l_i = 2r + 1 + (n -3)$ and $d \in \mathbb{N}$,
\begin{align*}
\left<\alpha_1,\ldots,\alpha_n\right>_{0,n,d} &\equiv
\int_{[\Mbar_{0,n}(X,d\ell)]^{virt}} e_1^*\alpha_1
\cdots e_n^*\alpha_n\\
&= (-1)^{(d - 1)(r + 1)} N_{l_1,\ldots,l_n} d^{n -
3}(\alpha_1.h^{r - l_1}) \cdots (\alpha_n.h^{r - l_n}).
\end{align*}
where $N_{l_1,\ldots,l_n}$ are recursively determined universal
constants. $N_{l_1,\ldots,l_n}$ are independent of $d$ and
$N_{l_1,\ldots,l_n} = 1$ for $n =2$ or $3$. All other (primary)
Gromov-Witten invariants with degree in $\mathbb{Z}\ell$ vanish.
\end{theorem}

This formula, together with some algebraic manipulations, implies
that for simple $\mathbb{P}^r$ flops the quantum corrections
attached to the extremal ray exactly remedy the defect caused by
the classical product for any $r \in \mathbb{N}$ and the big
quantum products restricted to exceptional curve classes  are
invariant under simple ordinary flops. Note that there are Novikov
variables $\q$ involved in these transformations (c.f.\ Remark
\ref{anal-cont}), and
\[
 \T (\q^\beta) = \q^{\T \beta}.
\]

The proof has two ingredients: Localization and the divisor
relations. Localization has been widely used in calculating
Gromov--Witten invariants. For genus zero one-pointed descendent
invariants twisted by a direct sum of negative line bundles, this
was carried out in \cite{LLY1} and \cite{Givental} in the context
of the study of {\it mirror symmetry}. The divisor relations
studied in \cite{LP} gives a {\it reconstruction theorem}, which
allows us to go from one-point invariants to multiple-point ones.

To achieve the invariance of big quantum product, non-extremal
curve classes need to be analyzed. The main purpose of \S\ref{s:4}
is to reduce the case of general $X$ to the local case. Briefly,
the degeneration formula expresses $\langle \alpha \rangle^X$ in
terms of relative invariants $\langle \alpha_1 \rangle ^{(Y, E)}$
and $\langle \alpha_2 \rangle ^{(\tilde E, E)}$, where $Y \to X$
is the blow-up of X over $Z$ and $\tilde{E} =
\mathbb{P}_Z(N_{Z/X}\oplus \mathscr{O})$. Similarly for $X'$, one
has $Y', \tilde{E}', E'$. By definition of ordinary flops, $Y =
Y'$ and $E = E'$. It is possible to match all output on the part
of $(Y,E)$ from $X$ and $X'$. Thus, the problem is transformed to
one for the relative cases of $(\tilde{E},E)$ and
$(\tilde{E}',E)$. Following ideas in the work of D.~Maulik and
R.~Pandharipande \cite{MP}, a further reduction from relative
invariants to absolute invariants is made. The problem is thus
reduced to
$$
X = \tilde E = \mathbb{P}_{\mathbb{P}^r}(\mathscr{O}(-1)^{\oplus(r
+ 1)}\oplus \mathscr{O}),
$$
which is a semi-Fano projective bundle.

\begin{remark} \label{r:1.8}
For simple flops, we may and will consider only cohomology
insertions of {\it real even degrees} throughout all our
discussions on GW invariants. This is allowed since $\tilde E$ has
only algebraic classes and any real odd degree insertion must go
to the $Y$ side after degeneration.
\end{remark}

The proof of the local case is carried out in \S\ref{s:5} by
exploring the compatibility of functional equations of $n$-point
functions under the reconstruction procedure of genus zero
invariants. It is easy to see that the Mori cone
$$
NE(X) = \mathbb{Z}_+\ell \oplus \mathbb{Z}_+\gamma
$$
with $\ell$ the line class in $Z$ and $\gamma$ the fiber line
class of $X = \tilde E \to Z$. The proof is based on an induction
on $d_2$ and $n$ with degree $\beta = d_1 \ell + d_2 \gamma$. The
case $d_2 = 0$ is handled by Theorem~\ref{t:1.7}. For $d_2
> 0$, the starting case, namely the one-point invariant, is again
based on localization technique on semi-Fano toric manifolds
\cite{Givental} and \cite{LLY3}.

\begin{theorem}[Functional equations for local models]
Consider an $n$-point function on $X =
\mathbb{P}_{\mathbb{P}^r}(\mathscr{O}(-1)^{\oplus(r + 1)}\oplus
\mathscr{O})$,
$$
\langle \alpha \rangle = \sum_{\beta \in NE(X)} \langle \alpha_1,
\ldots, \alpha_n \rangle_\beta\, \q^\beta
$$
where $\alpha_i$ lies in the span of cohomology classes in $X$ and
descendents of (push-forward of) cohomology classes in $E$. For
$\beta = d_1 \ell + d_2 \gamma$, the summands are non-trivial only
for a fixed $d_2$. If $d_2 \ne 0$ then
$$
\T \langle \alpha \rangle^X \cong \langle \T \alpha \rangle^{X'}.
$$
\end{theorem}

(Here $\cong$ stands for equality up to analytic continuations.)
Combining all the previous results Theorem~\ref{t:1.2} is proved.

\begin{remark}
Concerning ingredient (1), it is very important to understand the
{\it closure} of ordinary flops. To the authors' knowledge, no
serious attempt was made toward a higher dimensional version of
(1) except some much weaker topological results \cite{Wang3}. Even
in dimension three, the only known proof of (1) relies on the
minimal model theory and classifications of terminal
singularities. It is desirable to have a direct proof in the
symplectic category. Such a proof should shed important light
toward the higher dimensional cases. Our main theorem applies to
$K$-equivalent maps that are composite of simple ordinary flops
and their limits.
\end{remark}

As an application of the construction of ordinary flops in
\S\ref{s:1}, we discuss (twisted) Mukai flops in \S\ref{s:6}. Some
new understanding of the local geometry of Mukai flops is
presented and this leads to a proof of Theorem~\ref{t:1.3}.
Theorem~\ref{t:1.3} can also be interpreted as a generalization of
a local version of Huybrechts' results on hyper-K\"ahler manifolds
\cite{Huyb1}, with the flexibility of allowing the base $S$ to be
any smooth variety. As in the hyper-K\"ahler case, it also implies
that the correspondence induced by the fiber product
$$
[X\times_{\bar X} X'] = [\bar\Gamma_f] + [Z\times_S Z'] \in A^*(X
\times X')
$$
is the one which gives an isomorphism of Chow motives.

Besides dimension three \cite{LiRu} and the hyper-K\"ahler case
\cite{Huyb1}, our results provide the first known series of
examples in all high dimensions which support
Conjecture~\ref{c:1.4}.

\subsection{Acknowledgements}
We would like to thank A.~Givental, C.-H.~Liu, D.~Maulik, Y.~Ruan,
S.-T.~Yau and J.~Zhou for useful discussions. C.-L.~W.~is grateful
to C.-S.~Lin and J.~Yu for their encouragement.

We are grateful to the National Center for Theoretic Sciences
(NCTS, Taiwan) for providing stimulating and delightful
environment which makes the collaboration possible.

Y.-P.~L.~is partially supported by NSF and AMS Centennial
Fellowship. H.-W.~L.~is partially supported by NSC. C.-L.~W.~is
partially supported by NSC and the NCTS Chern Fellowship.

\section{Ordinary flops} \label{s:1}

\subsection{Ordinary $\mathbb{P}^r$ flops.}
Let $\psi: X \to \bar{X}$ be a flopping contraction as defined in
\S\ref{s:0.1}. Our first task is to show that the corresponding
algebraic ordinary flop $X \dashrightarrow X'$ exists. The
construction of the desired flop is rather straightforward. First
blow up $X$ along $Z$ to get $\phi : Y \rightarrow X$. The
exceptional divisor $E$ is a $\mathbb{P}^r \times
\mathbb{P}^{r}$-bundle over $S$. The key point is that one may
blow down $E$ along another fiber direction $\phi' : Y \rightarrow
X'$, with exceptional loci $\bar\psi' : Z' = \mathbb{P}_S(F')
\rightarrow S$ for $F'$ another rank $r + 1$ vector bundle over
$S$ and also $N_{Z'/X'}|_{\bar\psi'-{\rm fiber}} \cong
\mathscr{O}_{\mathbb{P}^{r}}(-1)^{\oplus(r + 1)}$. We start with
the following elementary lemma.

\begin{lemma}
Let $p:Z=\mathbb{P}_S(F)\to S$ be a projective bundle over $S$ and
$V\to Z$ a vector bundle such that $V|_{p^{-1}(s)}$ is trivial for
every $s\in S$. Then $V \cong p^*F'$ for some vector bundle $F'$
over $S$.
\end{lemma}

\Proof Recall that $H^i(\mathbb{P}^r,\mathscr{O})$ is zero for
$i\ne 0$ and $H^0(\mathbb{P}^r, \mathscr{O}) \cong \mathbb{C}$. By
the theorem on Cohomology and Base Change we conclude immediately
that $p_*\mathscr{O}(V)$ is locally free over $S$ of the same rank
as $V$. The natural map between locally free sheaves
$p^*p_*\mathscr{O}(V)\to \mathscr{O}(V)$ induces isomorphisms over
each fiber and hence by the Nakayama Lemma it is indeed an
isomorphism. The desired $F'$ is simply the vector bundle
associated to $p_*\mathscr{O}(V)$. \Endproof

Now apply the lemma to $V =
\mathscr{O}_{\mathbb{P}_S(F)}(1)\otimes N_{Z/X}$, and we conclude
that
$$
N_{Z/X}\cong \mathscr{O}_{\mathbb{P}_S(F)}(-1)\otimes
\bar\psi^*F'.
$$
Therefore, on the blow-up $\phi:Y={\rm Bl}_Z X \to X$,
\[
  N_{E/Y} = \mathscr{O}_{\mathbb{P}_Z(N_{Z/X})}(-1).
\]
From the Euler sequence which defines the universal sub-line
bundle we see easily that $\mathscr{O}_{\mathbb{P}_Z(L\otimes
F)}(-1) = \bar\phi^*L\otimes \mathscr{O}_{\mathbb{P}_Z(F)}(-1)$
for any line bundle $L$ over $Z$. Since the projectivization
functor commutes with pull-backs, we have
$$
E = \mathbb{P}_Z(N_{Z/X})\cong\mathbb{P}_Z(\bar\psi^*F') =
\bar\psi^*\mathbb{P}_S(F')=\mathbb{P}_S(F)\times_S\mathbb{P}_S(F').
$$

For future reference we denote the projection map
$Z':=\mathbb{P}_S(F')\to S$ by $\bar\psi'$ and $E\to Z'$ by
$\bar\phi'$. The various sets and maps are summarized in the
following commutative diagram. \small
$$
\xymatrix{&E = \mathbb{P}_S(F) \times_S \mathbb{P}_S(F')\subset Y
\ar[ld]_{\bar\phi}
\ar[rd]^{\bar\phi'}&\\
Z = \mathbb{P}_S(F)\subset X\;\ar[rd]_{\bar\psi}& &Z' =
\mathbb{P}_S(F')\ar[ld]^{\bar\psi'}\\
&S\subset \bar X}
$$ \normalsize
with normal bundle of $E$ in $Y$ being
\begin{align*}
N_{E/Y} &= \mathscr{O}_{\mathbb{P}_Z(N_{Z/X})}(-1) =
\mathscr{O}_{\mathbb{P}_Z(\mathscr{O}_Z(-1)\otimes
\bar\psi^*F')}(-1)
\\
&= \bar\phi^*\mathscr{O}_{\mathbb{P}_S(F)}(-1)\otimes
\mathscr{O}_{\mathbb{P}_Z(\bar\psi^*F')}(-1)\\
&= \bar\phi^*\mathscr{O}_{\mathbb{P}_S(F)}(-1)\otimes
\bar\phi'^*\mathscr{O}_{\mathbb{P}_S(F')}(-1).
\end{align*}

\begin{remark}
Notice that the bundles $F$ and $F'$ are uniquely determined up to
a twisting by a line bundle. Namely, the pair $(F, F')$ is
equivalent to $(F\otimes L, F'\otimes L^*)$ for any line bundle
$L$ on $S$.
\end{remark}

The next step is to show that there is a blow-down map $\phi':Y
\to X'$ which contracts the left ruling of $E$ and restricts to
the projection map $\bar\phi':E\to Z'$. The existence of the
contraction $\psi:X\to\bar X$ is essential here. Let us denote a
line in the left ruling by $C_Y$ such that $\phi(C_Y) = C$.

\begin{proposition} \label{existence}
Ordinary $\mathbb{P}^r$ flops exist.
\end{proposition}

\Proof Firstly, we will show that $C_Y$ is $K_Y$-negative. From
the exact sequence $0\to T_C\to T_X|_C \to N_{C/X}\to 0$ and
$N_{C/X}\cong \mathscr{O}_C(1)^{\oplus (r - 1)} \oplus
\mathscr{O}_C(-1)^{\oplus (r + 1)} \oplus \mathscr{O}_C^{\dim S}$,
we find that
$$
(K_X.C) = 2g(C) - 2 - ((r - 1) - (r + 1)) = 0.
$$
Together with $K_Y = \phi^*K_X + rE$, we get
$$
(K_Y.C_Y) = (K_X.C) + r(E.C_Y) = - r < 0.
$$

Next we will show $C_Y$ is extremal, i.e.~it has supporting (big
and nef) divisors. Let $H$ be a very ample divisor on $X$ and $L$
a supporting divisor for $C$ (e.g.\ take $L = \phi^*\bar H$ for an
ample divisor $\bar H$ on $\bar X$). Let $c = (H.C)$, then
$\phi^*H + cE$ has type $(0,-c)$ on each $\mathbb{P}^r\times
\mathbb{P}^r$ fiber of $E$. The divisor
$$
k\phi^*L - (\phi^*H + cE)
$$
is clearly big and nef for large $k$ and vanishes precisely on the
class $[C_Y]$. Thus $C_Y$ is a $K_Y$-negative extremal ray and the
contraction morphism $\phi': Y \to X'$ fits into\small
$$
\xymatrix{Y\ar[rr]^{\phi'}\ar[rd]_{\psi\circ\phi}&
& X'\ar[ld]^{\psi'}\\
& \bar X}
$$\normalsize
by the cone theorem on $Y \to \bar X$ (c.f.\ \cite{KoMo2}).
$X\dashrightarrow X'$ is then the desired flop. \Endproof

\begin{remark} \label{supp_div}
Notice that $(K_X.C) = 0$, $(K_{X'}.C') = 0$ ($C'$ is a line in
the fiber of $Z'\to S$) and $\phi^*K_X = \phi'^*K_{X'}$
($K$-equivalence).

It is clear from the proof that for the existence of $\phi'$ one
needs only the (weaker) assumption that $C$ is extremal instead of
the existence of the contraction $\psi:X \to \bar X$. However,
since $(K_X.C) = 0$ these two are indeed equivalent by the cone
theorem.
\end{remark}

\subsection{Local models} \label{s:1-2}
In general, without assuming the existence of $\psi$, (i) and (ii)
are not sufficient to construct $\phi'$ in the projective
category. This is well known already in the case of Atiyah flop
($r = 1$ and $S = \{\hbox{pt}\}$). In the analytic category
results of Cornalba \cite{Cornalba} do imply the contractibility
of $\psi$, $\phi'$ and $\psi'$ hence lead to the existence of
analytic ordinary $\mathbb{P}^r$ flops under (i) and (ii). The
situation is particularly simple in the case of {\it local models}
which we now describe.

Consider a complex manifold $S$ and two holomorphic vector bundles
$F\to S$ and $F'\to S$. Let $\bar\psi: Z:=\mathbb{P}_S(F)\to S$
and $\bar\psi': Z' :=\mathbb{P}_S(F')\to S$ be the induced
morphisms and let $E = \mathbb{P}_S(F)\times_S \mathbb{P}_S(F')$
with two projections $\bar\phi:E\to Z$ and $\bar\phi':E\to Z'$.
Let $Y$ be the total space of $N := \bar\phi^*
\mathscr{O}_Z(-1)\otimes \bar\phi'^* \mathscr{O}_{Z'}(-1)$ with
$E$ the zero section. It is clear that $N_{E/Y} = N$. There is a
contraction diagram \small
$$
\xymatrix{ &
E\;\ar[ld]_>>>>{\pi_1=\bar\phi}\ar[rd]|\hole^>>>>{\bar\phi'=\pi_2}
\ar@{^{(}->}[r]^j &
Y\;\ar[ld]_>>>>\phi\ar[rd]^>>>>{\phi'}\\
Z\;\ar[rd]_>>>>>>>>{\bar\psi}\ar@{^{(}->}[r]^i &
X\;\ar[rd]_>>>>>>>>>>{\psi} &
Z'\;\ar[ld]|\hole^>>>>>>>>>>{\bar\psi'} \ar@{^{(}->}[r]^{i'} &
X'\;\ar[ld]^>>>>>>>>{\psi'}\\
&S\;\ar@{^{(}->}[r]^{j'}& \bar X& }
$$ \normalsize
in the analytic category, with $X$ (resp.~$X'$) being the total
space of $\mathscr{O}_{\mathbb{P}_S(F)}(-1)\otimes \bar\psi^*F'$
(resp.~$\mathscr{O}_{\mathbb{P}_S(F')}(-1)\otimes \bar\psi'^*F$).

First of all, the discussion in \S1.1 implies that $\phi$ and
$\phi'$ are simply the blow-up maps along $Z$ and $Z'$
respectively. For $\psi$ and $\psi'$, when $S$ reduces to a point
the existence of contraction morphism $g: (Y, E) \to (\bar X, {\rm
pt})$ is a classical result of Grauert since $N_{E/Y}$ is a
negative line bundle. From the universal property the induced maps
$\psi$ and $\psi'$ are then analytic. For $S$ a small Stein open
set, $g: (Y, E) \to (\bar X, S)$, as well as $\psi$ and $\psi'$,
also exists since the whole picture is a trivial product with $S$.
The general case follows from patching the local data over an open
cover of $S$. In summary the local analytic model of an ordinary
$\mathbb{P}^r$ flop is a locally trivial family (over $S$) of
simple ordinary $\mathbb{P}^r$ flops.

It is convenient to consider compactified local models $\tilde X$,
$\tilde Y$ etc.~by adding the common infinity divisor $E_\infty
\cong E$ to $X$, $Y$ etc.~respectively. Denote by
$$
p: \tilde X = \mathbb{P}_Z(N_{Z/X} \oplus \mathscr{O}_Z) \to Z.
$$

\begin{proposition}
If $S$ is projective, for any bundles $F$, $F'$ of rank $r + 1$
the compactified local models of $\mathbb{P}^r$ flops exist in the
projective category.
\end{proposition}

\Proof $\tilde X$ is clearly projective and $E_\infty$ is
$p$-ample. By Remark \ref{supp_div} and Proposition
\ref{existence} we only need to construct a supporting divisor $L$
for the fiber line of $\bar \psi: Z \to S$. Let $H$ be ample in
$S$ then $\bar\psi^*H$ is a supporting divisor for the fiber line
in $Z$. Hence we may take $L := p^* \bar\psi^* H + E_\infty$.
\Endproof

The projective local models will be used extensively in
\S\ref{s:4}--\S\ref{s:6}.

\section{Correspondences and motives}  \label{s:2}

\subsection{Grothendieck's category of Chow motives}
General references of Chow motives can be found in \cite{Manin}
and \cite{Fulton}. Let $\M$ be the category of Chow motives (over
$\mathbb{C}$). For each smooth variety $X$, one associates an
object $\hat{X}$ in $\M$. The morphisms are given by
correspondences
$$
{\rm Hom}_{\M}(\hat X_1,\hat X_2) = A^*(X_1\times X_2).
$$
For $U\in A^*(X_1\times X_2)$, $V\in A^*(X_2\times X_3)$, let
$p_{ij}: X_1\times X_2\times X_3 \to X_i\times X_j$ be the
projection maps. The composition law is given by
$$
V\circ U = {p_{13}}_*(p_{12}^*U.p_{23}^*V).
$$

A correspondence $U$ has associated maps on Chow groups:
$$
U: A^*(X_1)\to A^*(X_2);\quad a\mapsto {p_2}_*(U.p_1^*a)
$$
as well as induced maps on $T$-valued points ${\rm Hom}(\hat
T,\hat X_i)$:
$$
U_T: A^*(T\times X_1) \mathop{\longrightarrow}^{U\circ}
A^*(T\times X_2).
$$
Then we have {\it Manin's identity principle}: Let $U,\, V\in {\rm
Hom}(\hat X, \hat X')$. Then $U = V$ if and only if $U_T = V_T$
for all $T$. (Since $U = U_X(\Delta_X) = V_X(\Delta_X) = V$.)

\begin{theorem}
For an ordinary $\mathbb{P}^r$ flop $f:X\dashrightarrow X'$, the
graph closure $\T := [\bar\Gamma_f]$ induces $\hat X \cong \hat
X'$ via $\T^*\circ \T = \Delta_X$ and $\T\circ \T^* =
\Delta_{X'}$.
\end{theorem}

\Proof For any $T$, ${\rm id}_T\times f: T\times X \dashrightarrow
T\times X'$ is also an ordinary $\mathbb{P}^r$ flop. Hence to
prove that $\T^*\circ \T = \Delta_X$, by the identity principle,
we only need to show that $\T^*\T = {\rm id}$ on $A^*(X)$ for any
ordinary $\mathbb{P}^r$ flop. From the definition of pull-back,
$$
\T W = p'_*(\bar\Gamma_f.p^*W) = \phi'_*\phi^*W.
$$
We also have the formulae for pull-back from the intersection
theory (c.f.~\cite{Fulton}, Theorem 6.7, Blow-up formula):
$$
\phi^*W = \tilde{W} + j_*\big(c(\mathscr{E}).\bar\phi^*s(W\cap
Z,W)\big)_{\dim W}
$$
where $\tilde{W}$ is the proper transform of $W$ in $Y$ and
$\mathscr{E}$ is the excess normal bundle defined by
\begin{equation} \label{e:excess}
 0\to N_{E/Y}\to \phi^* N_{Z/X} \to \mathscr{E}\to 0
\end{equation}
and $s(W \cap Z, W)$ is the relative Segre class. The key
observation is that the error term is lying over $W\cap Z$.

Let $W\in A_k(X)$. By Chow's moving lemma we may assume that $W$
intersects $Z$ transversally, so
$$
\ell := \dim W\cap Z = k + (r + s) - (r + r + s + 1) = k - r - 1.
$$
Since $\dim \phi^{-1}(W\cap Z) = \ell + r = k - 1 < k$, the error
term in the pull-back formula must be zero and we get $\phi^*W =
\tilde{W}$. Hence $\T W = W'$, the proper transform of $W$ in
$X'$. Notice that $W'$ is almost never transversal to $Z'$.

Let $B$ be an irreducible component of $W\cap Z$ and $\bar B =
\bar\psi(B) \subset S$ with dimension $\ell_B \le \ell$. Notice
that $W'\cap Z'$ has irreducible components $\{B' :=
\bar\psi'^{-1}(\bar B)\}_{B'}$ (different $B$ with the same $\bar
B$ will give rise to the same $B'$).

Let $\phi'^*W' = \tilde{W} + \sum E_{B'}$, where $E_{B'}$ varies
over irreducible components lying over $B'$, hence $E_{B'} \subset
\bar\phi'^{-1}\bar\psi'^{-1}(\bar B)$, a $\mathbb{P}^r\times
\mathbb{P}^{r}$ bundle over $\bar B$. For the generic point $s\in
\psi(\phi(E_{B'})) \subset \bar B$, we thus have
$$
\dim E_{B',s} \ge k - \ell_B = r + 1 + (\ell - \ell_B) > r.
$$
In particular, $E_{B',s}$ contains positive dimensional fibers of
$\phi$ (as well as $\phi'$). Hence $\phi_*(E_{B'}) = 0$ and
$\T^*\T W = W$.

By the same argument we have also that $\T\circ \T^* =
\Delta_{X'}$, thus the proof is completed. \Endproof

\begin{remark}
For a general ground field $k$, if the flop diagram under
consideration is defined over $k$ then the theorem works for
motives over $k$.
\end{remark}

\begin{corollary}\label{poincare}
Let $f:X \dashrightarrow X'$ be a $\mathbb{P}^r$ flop. If $\dim
\alpha_1 + \dim \alpha_2 = \dim X$, then
$$
(\T \alpha_1.\T \alpha_2) = (\alpha_1.\alpha_2).
$$
That is, $\T$ is an isometry with respect to $(-.-)$.
\end{corollary}

\Proof We may assume that $\alpha_1$, $\alpha_2$ are transversal
to $Z$. Then
\begin{align*}
(\alpha_1.\alpha_2) &= (\phi^*\alpha_1.\phi^*\alpha_2) =
((\phi'^*\T\alpha_1 - \xi).\phi^*\alpha_2)\\ &=
((\phi'^*\T\alpha_1).\phi^*\alpha_2) = (\T\alpha_1.(\phi'_*
\phi^*\alpha_2)) = (\T\alpha_1.\T\alpha_2).
\end{align*}
Here we use the fact proved in the above theorem that $\xi$ has
positive fiber dimension in the $\phi$ direction. \Endproof

Thus for ordinary flops, $\T^{-1} = \T^*$ both in the sense of
correspondences and Poincar\'e pairing.

\begin{remark}
It is an easy fact that if $X =_K X'$ then $X$ and $X'$ are
isomorphic in codimension one and in particular the graph closure
gives canonical isomorphisms $\T$ on $A^1(X) \cong A^1(X')$ and
$A_1(X) \cong A_1(X')$ respectively. In this more general setting,
the above proof still implies that the Poincar\'e pairing on
$A^1\times A_1$ (and $H^2 \times H_2$) is preserved under $\T$.
\end{remark}

\subsection{Triple product for simple flops}
Let $f: X\dashrightarrow X'$ be a simple $\mathbb{P}^r$ flop with
$S$ being a point. Let $h$ be the hyperplane class of $Z =
\mathbb{P}^r$ and $h'$ be the hyperplane class of $Z'$. Let also
$x = \bar\phi^*h = [h\times \mathbb{P}^r]$, $y = \bar\phi'^*h' =
[\mathbb{P}^r\times h']$ in $E = \mathbb{P}^r\times \mathbb{P}^r$.

\begin{lemma}\label{linear-class}
For classes inside $Z$, we have
$$
\phi^*[h^l] = j_*(x^ly^r - x^{l+1}y^{r - 1} + \cdots + (-1)^{r -
l}x^ry^l).
$$
Hence by symmetry we get $\T[h^l] = (-1)^{r - l}[h'^l]$. In
particular, $\T[C] = -[C']$.
\end{lemma}

\Proof Recall that
$$
N_{E/Y} = \mathscr{O}_{\mathbb{P}^r\times\mathbb{P}^r}(-1,-1) :=
\bar\phi^*\mathscr{O}_{\mathbb{P}^r}(-1)\otimes
\bar\phi'^*\mathscr{O}_{\mathbb{P}^r}(-1)
$$
and $N_{Z/X} = \mathscr{O}_{\mathbb{P}^r}(-1)^{\oplus (r + 1)}$.
From \eqref{e:excess},
$$
c(\mathscr{E}) = (1 - x)^{r + 1}(1 - x - y)^{-1}.
$$
Taking degree $r$ terms from both sides, we have
\begin{align*}
c_r(\mathscr{E}) &= [(1 - x)^{r + 1}(1 - (x + y))^{-1}]_{(r)}\\
&= (x + y)^r - C^{r + 1}_1(x + y)^{r - 1}x + \cdots +
(-1)^r C^{r + 1}_r x^r\\
&= (x + y)^{-1}((x + y) - x)^{r + 1} - (-1)^{r + 1}x^{r + 1})\\
&= (y^{r + 1} - (-1)^{r + 1}x^{r + 1})/(y + x)\\
&= y^r - y^{r-1}x  + y^{r-2}x^2 - \cdots + (-1)^rx^r.
\end{align*}
The basic pull-back formula (\cite{Fulton}, Proposition 6.7) then
implies that
$$
\phi^*[h^l] = j_*(c_r(\mathscr{E}).\bar\phi^*[h^l]) =
j_*(c_r(\mathscr{E}).x^l) = j_*\sum\nolimits_{t = 0}^r (-1)^t y^{r
- t}x^{t + l}.
$$
If $t + l \geq r + 1$ then $y^{r - t}x^{t + l} = 0$. The result
follows. \Endproof

\begin{lemma}
For a class $\alpha\in H^{2l}(X)$ with $l \le r$, let $\alpha' =
\T \alpha$ in $X'$. Then
$$
\phi'^*\alpha' = \phi^*\alpha + (\alpha.h^{r - l})\,j_*\frac{x^l -
(-y)^l}{x + y}.
$$
\end{lemma}

\Proof Since the difference $\phi'^*\alpha' - \phi^*\alpha$ has
support in $E$, we may write
$$
\phi'^*\alpha' = \phi^*\alpha + j_*(a_1x^{l - 1} + \cdots +
a_kx^{l - k}y^{k - 1} + \cdots + a_ly^{l -1}).
$$
By intersecting this equation with $x^{r - l}y^r$ in $X$ and
noticing that $E\sim -(x + y)$ on $E$, we get by the projection
formula
$$
0 = \phi^*\alpha.x^{r - l}y^r - a_1x^{l - 1}(x + y)x^{r - l}y^r =
(\alpha.h^{r - l}) - a_1.
$$
Similarly by intersecting with $x^{r - l + 1}y^{r - 1}$ we get
$$
0 = - a_1x^{l - 1}(x + y)x^{r - l + 1}y^{r - 1} - a_2x^{l - 2}(x +
y)x^{r - l + 1}y^{r - 1} = -a_1 - a_2.
$$

Continuing in this way by intersecting with $x^py^q$ with $p + q =
2r - l$ we get $a_k = (-1)^{k - 1}(\alpha.h^{r - l})$ for all $k =
1,\ldots, l$. This proves the lemma. \Endproof

These formulae allow us to compare the triple products of classes
in $X$ and $X'$:

\begin{proposition}\label{classical}
For a simple $\mathbb{P}^r$-flop $f:X \dashrightarrow X'$, let
$\alpha_i \in H^{2 l_i}(X)$, with $l_i \le r$, $l_1 + l_2 + l_3 =
\dim X = 2r + 1$. Then
$$
(\T\alpha_1.\T\alpha_2.\T\alpha_3) = (\alpha_1.\alpha_2.\alpha_3)
+ (-1)^r (\alpha_1.h^{r -l_1})(\alpha_2.h^{r -l_2})(\alpha_3.h^{r
-l_3}).
$$
\end{proposition}

\Proof The proof consists of straightforward computations.
\begin{align*}
&(\T\alpha_1.\T\alpha_2.\T\alpha_3) =
(\phi'^*\T\alpha_1.\phi'^*\T\alpha_2.\phi'^*\T\alpha_3)\\
&= \left(\phi^*\alpha_1 + (\alpha_1.h^{r -l_1})j_*\frac{x^{l_1} -
(-y)^{l_1}}{x + y}\right) \left(\phi^*\alpha_2 + (\alpha_2.h^{r
-l_2})
 j_*\frac{x^{l_2} - (-y)^{l_2}}{x + y}\right)\\
&\qquad \qquad \qquad \times \left(\phi^*\alpha_3 + (\alpha_3.h^{r
- l_3})j_*\frac{x^{l_3} - (-y)^{l_3}}{x + y}\right).
\end{align*}

Among the resulting eight terms, the first term is clearly equal
to $\alpha_1.\alpha_2.\alpha_3$.

For those three terms with two pull-backs like
$\phi^*\alpha_1.\phi^*\alpha_2$, the intersection values are zero
since the remaining part necessarily contains the $\phi$ fiber
(from the formula the power in $y$ is at most $l_3 - 1$).

The term with $\phi^*\alpha_1$ and two exceptional parts
contributes
\begin{align*}
&\phi^*\alpha_1.j_*\frac{x^{l_2} -
(-y)^{l_2}}{x + y}.j_*\frac{x^{l_3} - (-y)^{l_3}}{x + y}\\
&= -\phi^*\alpha_1.j_* \big((x^{l_2} -(-y)^{l_2})(x^{l_3 - 1}
+x^{l_3 - 2}(-y) + \cdots + (-y)^{l_3 - 1})\big)
\end{align*}
times $(\alpha_2.h^{r - l_2})(\alpha_3.h^{r - l_3})$. The terms
with non-trivial contribution must contain $y^r$, hence there is
only one such term, namely (notice that $l_1 + l_2 +l_3 = 2r + 1$)
$$
-(-y)^{l_2}\times x^{l_3 - 1 -(r - l_2)}(-y)^{r - l_2} = -(-1)^r
x^{r-l_1}y^r
$$
and the contribution is $(-1)^r(\alpha_1.h^{r
-l_1})(\alpha_2.h^{r-l_2}) (\alpha_3.h^{r -l_3})$. There are three
such terms.

It remains to consider the term of triple product of three
exceptional parts. It is $(\alpha_1.h^{r -l_1})(\alpha_2.h^{r
-l_2}) (\alpha_3.h^{r -l_3})$ times
$$
(x^{l_1} - (-y)^{l_1})(x^{l_2} - (-y)^{l_2})(x^{l_3 - 1} + x^{l_3-
2}(-y) + \cdots +(-y)^{l_3 - 1}).
$$
The terms with non-trivial values are precisely multiples of $x^r
y^r$. Since $l_1 + l_2 > r$, there are two such terms
$$
-x^{l_1} (-y)^{l_2} \times x^{r -l_1}(-y)^{l_3- 1 -(r -l_1)} -
x^{l_2} (-y)^{l_1} \times x^{r -l_2}(-y)^{l_3 - 1 -(r - l_2)}
$$
which give $-2(-1)^r$. Summing together we then finish the proof.
\Endproof

\subsection{Motives and ordinary flips}
Results in \S1 and \S2 extend straightforwardly to the case of
ordinary flips. Before we move to quantum corrections for ordinary
flops, we shall summarized here the classical aspects, especially
the motivic aspects, of ordinary flips. The proofs are identical
with the flop case and are thus omitted.

Consider $(\psi, \bar \psi): (X, Z) \to (\bar X, S)$ a
log-extremal contraction as before. $\psi$ is an ordinary $(r,r')$
flipping contraction if
\begin{itemize}
\item[(i)] $Z = \mathbb{P}_S(F)$ for some rank $r + 1$ vector
bundle $F$ over $S$,

\item[(ii)] $N_{Z/X}|_{Z_s} \cong
\mathscr{O}_{\mathbb{P}^r}(-1)^{\oplus(r' + 1)}$ for each
$\bar\psi$-fiber $Z_s$, $s\in S$.
\end{itemize}
Then the $(r, r')$ flip $f: X \dashrightarrow X'$ exists with
explicit local model as in \S\ref{s:1-2}.

In terms of the $K$-partial order within a birational class, $X
\le_K X'$ if and only if $r \le r'$. For $f$ a $(r, r')$ flip with
$r \le r'$, the graph closure $\T = [\bar\Gamma_f] \in A^*(X
\times X')$ identifies the Chow motive $\hat X$ of $X$ as a
sub-motive of $\hat X'$ which preserves also the Poincar\'e
pairing on cohomology groups.

More precisely, a self correspondence $p \in A^*(X\times X)$ is a
projector if $p^2 = p$. There is a natural {\it pseduo-abelian
extension} $\tilde\M$ of $\M$ to include all pairs $(X, p)$ as its
objects. $(X, p)$ is regarded as the {\it image} of $p$. Moreover,
$\hat X = (X, p)\oplus (X, 1 - p)$ in $\tilde\M$. With this
notion, for an ordinary $(r, r')$ flip $f:X\dashrightarrow X'$
with $r\le r'$, the graph closure $\T := [\bar\Gamma_f]$ induces
$\hat X \cong (X',p')$ via $\T^*\circ \T = \Delta_X$, where $p' =
\T\circ \T^*$ is a projector.

Since every {\it geometric cohomology theory} (a graded ring
functor $H^*$ with Poincar\'e duality, K\"unneth formula and a
cycle map $A^* \to H^*$ etc.) factors through $\tilde \M$, the
result also holds on such a specialized theory.

For simple $(r, r')$ flips (i.e.\ $S = {\rm pt}$) with $l \le
\min\{r, r'\}$,
$$
\phi^*[h^{r - l}] = j_*(x^{r - l}y^{r'} - x^{r - l + 1}y^{r' - 1}
+ \cdots + (-1)^{l}x^ry^{r' - l}).
$$
In particular $\T[h^{r - l}] = (-1)^{l}[h'^{r' - l}]$. For
$\alpha\in A^l(X)$ with $l \le \min\{r, r'\}$,
$$
\phi'^* \T\alpha = \phi^*\alpha + (\alpha.h^{r - l})\,j_*\frac{x^l
- (-y)^l}{x + y}.
$$

Let $\alpha_i\in H^{2l_i}(X)$, $1 \le i \le 3$ with $l_i \le
\min\{r, r'\}$, $l_1 + l_2 + l_3 = \dim X = r + r' + 1$. The
defect of the triple product is again given by
$$
(\T\alpha_1.\T\alpha_2.\T\alpha_3) = (\alpha_1.\alpha_2.\alpha_3)
+ (-1)^{r'} (\alpha_1.h^{r - l_1})(\alpha_2.h^{r -
l_2})(\alpha_3.h^{r - l_3}).
$$

\section{Quantum corrections attached to extremal rays} \label{s:3}
Proposition~\ref{classical} on triple products suggests that one
needs to correct the product structure by some contributions from
the extremal ray. In this section we show that for simple ordinary
flops the quantum corrections attached to the extremal ray exactly
remedy the defect of the ordinary product.

\subsection{Quantum cohomology} \label{s:3.1}
We use \cite{CoKa} as our general reference on moduli spaces of
stable maps, Gromov-Witten theory and quantum cohomology.

Let $\beta \in NE(X)$, the Mori cone of numerical classes of
effective one cycles. Let $\Mbar_{g,n}(X, \beta)$ be the moduli
space of $n$-pointed stable maps $f:(C; x_1, \ldots, x_n) \to X$
from a nodal cure $C$ with arithmetic genus $g(C) = g$ and with
degree $[f(C)] = \beta$. Let $e_i:\Mbar_{g,n}(X, \beta) \to X$ be
the evaluation morphism $f \mapsto f(x_i)$. The Gromov-Witten
invariant for classes $\alpha_i \in H^{*}(X)$, $1\le i \le n$, is
given by
$$
\left<\alpha_1,\ldots,\alpha_n\right>_{g,n,\beta} := \int_{[\bar
M_{g,n}(X,\beta)]^{virt}} e_1^*\alpha_1 \cdots e_n^*\alpha_n.
$$

The genus zero three-point functions (as formal power series)
$$
\left<\alpha_1,\alpha_2,\alpha_3\right> := \sum\nolimits_{\beta\in
A_1(X)}\left<\alpha_1, \alpha_2, \alpha_3\right>_{0,3,\beta}
\q^\beta
$$
together with the Poincar\'e pairing $(-,-)$ determine the {\it
small} quantum product.

More precisely, let $T = \sum t_iT_i$ with $\{T_i\}$ a cohomology
basis and $t_i$ being formal variables. Let $\{T^i\}$ be the dual
basis with $(T^i,T_j) = \delta_{ij}$. The (genus zero)
pre-potential combines all $n$-point functions together:
$$
\Phi(T) = \sum\nolimits_{n = 0}^\infty \sum\nolimits_{\beta \in
NE(X)} \frac{1}{n!} \left<T^n\right>_{\beta}\,\q^\beta,
$$
where $\big<T^n\big>_{\beta} = \left<T, \ldots,
T\right>_{0,n,\beta}$. The \emph{big quantum product} is defined by
$$
T_i *_t T_j = \sum\nolimits_k \Phi_{ijk} T^k
$$
where
$$
\Phi_{ijk} = \frac{\partial^3\Phi}{\partial t_i\partial t_j
\partial t_k} = \sum\nolimits _{n = 0}^\infty
\sum\nolimits_{\beta \in NE(X)} \frac{1}{n!} \left<T_i, T_j, T_k,
T^n\right>_\beta \q^\beta.
$$
The \emph{small quantum product} is defined to be the restriction
of $*_t$ to $t=0$ (the $n = 0$ part $\Phi_{ijk}(0)$).

In general, it is difficult to calculate $\Phi(T)$ and the big
quantum ring directly. When $X$ admits symmetries and $H^*(X)$ is
\emph{generated by divisors}, it is usually possible to use
localization techniques to calculate one point invariants with
gravitational descendents, or its generating function, the
$J$-function, defined as follows.
\begin{equation} \label{e:J}
 \begin{split}
 J_X(\q,z^{-1}) :=
 &\sum_{\beta \in NE(X)} \q^{\beta} J_X(\beta, z^{-1})
 \in H^*(X)[\![z^{-1}]\!][\![\q]\!] \\
 := &\sum_{\beta \in NE(X)} \q^{\beta}
 e_{1*}^X \left( \frac{1}{z(z-\psi)} \cap
 [\Mbar_{0,1}(X,\beta)]^{virt} \right)
 \end{split}
\end{equation}
Furthermore, the reconstruction theorem in \cite{LP} (also
\cite{BK}) implies that $J$-function actually determines the
entire generation function $\Phi(T)$.

\subsection{Analytic continuation}
Let $f:X\dashrightarrow X'$ be a simple $\mathbb{P}^r$ flop. Since
$X$ and $X'$ have the same Poincar\'e pairing under $\T$, in order
to compare their quantum products we only need to compare their
$n$-point functions. For three-point functions, write
$$
\left<\alpha_1,\alpha_2,\alpha_3\right> =
(\alpha_1.\alpha_2.\alpha_3) + \sum\nolimits_{d\in\mathbb{N}}
\left<\alpha_1, \alpha_2, \alpha_3\right>_{d\ell} \q^{d\ell} +
\sum\nolimits_{\beta \not\in \mathbb{Z}\ell}\left<\alpha_1,
\alpha_2, \alpha_3\right>_{\beta} \q^\beta.
$$
The difference $(\T\alpha_1.\T\alpha_2.\T\alpha_3) -
(\alpha_1.\alpha_2.\alpha_3)$ is already determined in last
section. The next step is to compute the middle term, namely
quantum corrections coming from the extremal ray $\ell = [C]$. The
third term will be discussed in later sections.

The virtual dimension of $\Mbar_{g,n}(X,d\ell)$ is given by
$$
(c_1(X).d\ell) + (2r + 1)(1 - g) + (3g - 3) + n.
$$
Since $(K_X.\ell) = 0$, for $g = 0$ we need only consider classes
$\alpha_i\in A^{l_i}(X)$ with $\sum_{i = 1}^n l_i = 2r + 1 + (n -
3)$. For $n = 3$ this is $2r + 1 = \dim X$.

\begin{theorem}\label{main}
For all $\alpha_i\in H^{2l_i}(X)$ with $1\le l_i \le r$, $\sum_{i
= 1}^n l_i = 2r + 1 + (n -3)$ and $d \in \mathbb{N}$,
\begin{align*}
\left<\alpha_1,\ldots,\alpha_n\right>_{0,n,d} &\equiv
\int_{[\Mbar_{0,n}(X,d\ell)]^{virt}} e_1^*\alpha_1
\cdots e_n^*\alpha_n\\
&= (-1)^{(d - 1)(r + 1)} N_{l_1,\ldots,l_n} d^{n -
3}(\alpha_1.h^{r - l_1}) \cdots (\alpha_n.h^{r - l_n}).
\end{align*}
where $N_{l_1,\ldots,l_n}$ are recursively determined universal
constants. $N_{l_1,\ldots,l_n}$ are independent of $d$ and
$N_{l_1,\ldots,l_n} = 1$ for $n =2$ or $3$. All other (primary)
Gromov-Witten invariants with degree in $\mathbb{Z}\ell$ vanish.
\end{theorem}

\begin{corollary}\label{q-inv}
Both the small and big quantum products restricted to exceptional
curve classes are invariant under simple ordinary flops. In fact
the three-point functions attached to the extremal ray exactly
remedy the defect caused by the classical product.
\end{corollary}

\Proof Since $(\T\alpha_i.h'^{(r - l_i)}) =
(-1)^{l_i}(\T\alpha_i.\T h^{r - l_i}) = (-1)^{l_i}(\alpha_i.h^{r -
l_i})$, for three point functions we get
\begin{align*}
&\left<\T\alpha_1,\T\alpha_2,\T\alpha_3\right> -
\left<\alpha_1,\alpha_2,\alpha_3\right> = (-1)^r (\alpha_1.h^{r -
l_1})(\alpha_2.h^{r - l_2})(\alpha_3.h^{r - l_3})\\
&\qquad + (\alpha_1.h^{r - l_1})(\alpha_2.h^{r -
l_2})(\alpha_3.h^{r - l_3})\left(\frac{(-1)^{2r + 1} \q^{\ell'}}{1
- (-1)^{r + 1} \q^{\ell'}} - \frac{\q^\ell}{1 - (-1)^{r + 1}
\q^\ell}\right).
\end{align*}

Under the correspondence $\T$, we shall identify $\q^{\ell'}$ with
$\q^{-\ell}$. Plug in this into the last bracket we get $1$ when
$r$ is odd and get $-1$ when $r$ is even. In both cases the right
hand side cancels out and then
$\left<\T\alpha_1,\T\alpha_2,\T\alpha_3\right> =
\left<\alpha_1,\alpha_2,\alpha_3\right>$. This proves the
statement on small quantum product.

For general $n = 3 + k$ point invariants with $k \ge 1$, we get
\begin{align*}
\left<\alpha_1,\ldots,\alpha_n\right> &=
N_{l_1,\ldots,l_n}(\alpha_1.h^{r - l_1})\cdots (\alpha_n.h^{r -
l_n})\sum_{d = 1}^\infty (-1)^{(d -
1)(r + 1)}d^k \q^{d\ell}\\
&= N_{l_1,\ldots,l_n}(\alpha_1.h^{r - l_1})\cdots (\alpha_n.h^{r -
l_n}) \left(\q^\ell \frac{d}{d\q^\ell} \right)^k \frac{(-1)^{r +
1}}{1 - (-1)^{r + 1}\q^\ell}.
\end{align*}
Similarly, since $(-1)^{\sum l_i} = (-1)^{k + 1}$,
$\left<\T\alpha_1,\ldots,\T\alpha_n\right>$ equals
$$
(-1)^{k + 1}N_{l_1,\ldots,l_n}(\alpha_1.h^{r - l_1})\cdots
(\alpha_n.h^{r - l_n}) \left(\q^{\ell'} \frac{d}{d\q^{\ell'}}
\right)^k \frac{(-1)^{r + 1}}{1 - (-1)^{r + 1}\q^{\ell'}}.
$$
Taking into account of
$$
\q^{-\ell} \frac{d}{d\q^{-\ell}} = - \q^\ell
\frac{d}{d\q^\ell}\quad \mbox{and}\quad \frac{1}{1 - (-1)^{r +
1}\q^{-\ell}} = 1 - \frac{1}{1 - (-1)^{r + 1}\q^\ell}
$$
we get $\left<\T\alpha_1,\ldots,\T\alpha_n\right> =
\left<\alpha_1,\ldots,\alpha_n\right>$ for all $k \ge 1$ ($n \ge
4$). The proof for the statement on big quantum product is thus
completed. \Endproof

To put the result into perspective, we interpret the change of
variable $\ell'$ by $-\ell$ in terms of analytic continuation over
the extended complexified K\"ahler moduli space.

Without lose of generality we illustrate this by writing out the
small quantum part. This is simply a word by word adoption of the
treatment in the $r = 1$ case (cf.\ \cite{Witten} \S5.5,
\cite{Morrison} \S4).

The quantum cohomology is parameterized by the complexified
K\"ahler class $\omega = B + iH$ with $\q^\beta = \exp(2\pi i
(\omega.\beta))$, where $B \in H^{1,1}_\mathbb{R} (X)$ and $H\in
\mathcal{K}_X$, the K\"ahler cone of $X$. For a simple
$\mathbb{P}^r$ flop $X \dashrightarrow X'$, $\T$ identifies
$H^{1,1}$, $A_1$ and the Poincar\'e pairing $(-,-)$ on $X$ and
$X'$. Then $\left<\alpha_1,\alpha_2,\alpha_3\right>^X$ restricted
to $\mathbb{Z}\ell$ converges in the region
$$
H^{1,1}_+ = \{\omega\,|\, (H.\ell) > 0\} \supset
H^{1,1}_\mathbb{R} \times i\,\mathcal{K}_X
$$
and equals
$$
(\alpha_1.\alpha_2.\alpha_3) + (\alpha_1.h^{r -
l_1})(\alpha_2.h^{r - l_2})(\alpha_3.h^{r - l_3})\frac{e^{2\pi
i(\omega.\ell)}}{1 - (-1)^{r + 1} e^{2\pi i(\omega.\ell)}}.
$$
This is a well-defined analytic function of $\omega$ on the whole
$H^{1,1}$, which defines the analytic continuation of
$\left<\alpha_1,\alpha_2,\alpha_3\right>^X$ from
$H^{1,1}_\mathbb{R} \times i\,\mathcal{K}_X$ to $H^{1,1}$.

Similarly, $\left<\T\alpha_1,\T\alpha_2,\T\alpha_3\right>^{X'}$
restricted to $\mathbb{Z}\ell'$ converges in the region
$$
\{\omega\,|\, (H.\ell') > 0\} = \{\omega\,|\, (H.\ell) < 0\} =
H^{1,1}_- \supset H^{1,1}_\mathbb{R} \times i\,\mathcal{K}_{X'}
$$
and equals
$$
(\T\alpha_1.\T\alpha_2.\T\alpha_3) - (\alpha_1.h^{r -
l_1})(\alpha_2.h^{r - l_2})(\alpha_3.h^{r - l_3})\frac{e^{-2\pi
i(\omega.\ell)}}{1 - (-1)^{r + 1} e^{-2\pi i(\omega.\ell)}}
$$
which is the analytic continuation of the previous one from
$H^{1,1}_+$ to $H^{1,1}_-$.\medskip

We introduce the notation $A \cong B$ for the two series $A$ and
$B$ when they can be analytically continued to each other.

\begin{remark}\label{anal-cont}
It was conjectured that the total series $\Phi^X_{ijk}$ converges
for $B\in \mathcal{K}_X$, at least for $B$ large enough, hence the
large radius limit goes back to the classical cubic product. The
Novikov variables $\{q^\beta\}_{\beta \in NE(X)}$ are introduced
to avoid the convergence issue.

Since $\mathcal{K}_X \cap \mathcal{K}_{X'} = \emptyset$ for
non-isomorphic $K$-equivalent models, the collection of K\"ahler
cones among them form a chamber structure. The conjectural
canonical isomorphism
$$
\T: H^*(X) \cong H^*(X')
$$
assigns to each model $X$ a {\it coordinate system} $H^*(X)$ of
the fixed $H^*$ and $\T$ serves as the {\it (linear) transition
function}. The conjecture asserts that $\Phi^{X}_{ijk}$ can be
{\it analytically continued} from $\mathcal{K}_{X}$ to
$\mathcal{K}_{X'}$ and agrees with $\Phi^{X'}_{ijk}$.
Equivalently, $\Phi_{ijk}$ is well-defined on $\mathcal{K}_X \cup
\mathcal{K}_{X'}$ which verifies the {\it functional equation}
$$
\T \Phi_{ijk}(\omega, T) \cong \Phi_{ijk}(\omega, \T T).
$$

For simple ordinary flops, this is verified from \S3 to \S5 for
each given cohomology insertions. The convergence has just been
verified for extremal rays and will be verified for local models
in \S5.
\end{remark}

\subsection{One-point functions with descendents}

In order to prove Theorem \ref{main}, we first reduce the problem
to one for projective spaces. Let
\[
 U_d := R^1\ft_* e_{n + 1}^* N
\]
be the obstruction bundle, where $N = N_{Z/X}$ and $\ft$ is the
forgetting morphism in
$$
\xymatrix{\Mbar_{0,n + 1}(\mathbb{P}^r,d)\ar[r]^{\qquad
e_{n + 1}}\ar[d]^{\ft} & \mathbb{P}^r\\
\Mbar_{0,n}(\mathbb{P}^r,d)}.
$$
It is well known (see e.g.\ \cite{CoKa}) that
\begin{equation} \label{e:virt}
[\Mbar_{0,n}(X, d\ell)]^{virt} =  e(U_d) \cap
[\Mbar_{0,n}(\mathbb{P}^r,d\ell)].
\end{equation}
Since $U_d$ is functorial under $\ft^*$, we use the same notation
for all $n$.

As explained earlier, we will start with the calculation of the
$J$-function \eqref{e:J}.
In our case
$$
J_X(d\ell,z^{-1}) \equiv e_{1*}^{\mathbb{P}^r}\frac{e(U_d)}{z(z -
\psi)}
$$
has been calculated: Let $\displaystyle{\Pd := (-1)^{(d - 1)(r +
1)} \frac{1}{(h + dz)^{r+1}}}$.

\begin{lemma}[\cite{LLY1}, also \cite{Givental}]
\label{l:J}
$$
J_X(d\ell,z^{-1}) = P_d.
$$
\end{lemma}

\begin{remark}
This calculation can be interpreted as quantum Lefschetz
hyperplane theorem for concave bundles over $\mathbb{P}^r$. From
this viewpoint, the ``mirror transformation'' from
$J_X(d\ell,z^{-1})$ to $\Pd$ is not needed since the rank of the
bundle $\mathscr{O}(-1)^{r+1}$ is greater than one. See
e.g.~\cite{ypL}.
\end{remark}

\begin{corollary}\label{one-point}
For $l + k = 2r - 1$, $1 \le l \le r$,
$$
\left<\tau_k h^l\right>_d = \frac{(-1)^{d(r + 1) + k}}{d^{k +
2}}\, C^{k + 1}_r
$$
where $C^k_r = {k!}/{r! (k-r)!}$. The invariant is zero if $l + k
\ne 2r - 1$ by dimensional constraints.
\end{corollary}

\Proof We start with
\[
A := \int_{\mathbb{P}^r} h^l. \Pd = \sum\nolimits_{k \ge 0}
\frac{1}{z^{k + 2}} \left<\tau_k h^l\right>_d.
\]
By Lemma~\ref{l:J}
$$
A = (-1)^{(d - 1)(r + 1)}\int_{\mathbb{P}^r} \frac{h^l}{(h +
dz)^{r + 1}} = \int_{\mathbb{P}^r} \frac{h^l}{d^{r + 1}z^{r +
1}}\left(1 + \frac{h}{dz}\right)^{-(r + 1)}.
$$
The result follows from the Taylor expansion and the elementary
fact that
$$
C^{-(r + 1)}_{r - l} =
(-1)^{k + (r + 1)}C^{k + 1}_r.
$$
\Endproof

\subsection{Multiple-point functions via divisor relations}

We recall the following rational equivalence in $A_*(\Mbar_{0,
n}(X, \beta)) \otimes \mathbb{Q}$ from \cite{LP}, Corollary 1: For
$L \in {\rm Pic}(X)$ and $i \ne j$,
\begin{equation}  \label{e:rec1}
\begin{split}
&e_i^*L\cap [\Mbar_{0,n}(X, \beta)]^{virt} \\
= &(e_j^*L + (\beta, L)\psi_j)\cap [\Mbar_{0,n}(X, \beta)]^{virt}
- \sum_{\beta_1 + \beta_2 =
\beta}(\beta_1,L)[D_{i,\beta_1|j,\beta_2}]^{virt},
\end{split}
\end{equation}
\begin{equation} \label{e:rec2}
\psi_i + \psi_j = [D_{i|j}]^{virt},
\end{equation}
where $[D_{i,\beta_1|j,\beta_2}]^{virt} \in A_*(\Mbar_{0, n}(X,
\beta))$ is the push-forward of the virtual classes of the
corresponding boundary divisor components
$$
D_{i,\beta_1|j, \beta_2} = \sum_{i \in A, j\in B;\ A \coprod B =
\{1,\dots, n\}} D(A, B; \beta_1, \beta_2)
$$
and
$$
D_{i|j} = \sum_{\beta_1 + \beta_2 = \beta} D_{i,\beta_1|
j,\beta_2}.
$$

Here is a simple observation which will be repeatedly used in the
sequel:

\begin{lemma}[Vanishing lemma]
Let $\mathbb{P}^r \subset X$ with $N_{\mathbb{P}^r/X} = \oplus_j
\mathscr{O}(-m_j)$, $ m_j \in \mathbb{N}$. Let $\ell$ be the line
class in $\mathbb{P}^r$. Then for $\deg T > r$ and $d \ne 0$,
$\langle \ldots, T \rangle_{d\ell} = 0$.
\end{lemma}

\Proof Since $[\Mbar_{0,n}(X, d\ell)]^{virt}$ equals $[\Mbar_{0,
n}(\mathbb{P}^r, d)]$ cut out by $e(U_{d})$, the evaluation
morphisms factor through $\mathbb{P}^r$. But then $e_n^* (T|_Z) =
0$. \Endproof

Here $\deg T := l$ if $T\in H^{2l}(X)$. As we had mentioned in the
introduction, only real even degree classes will be relevant
throughout our discussions.

\begin{proposition} \label{main-2}
For $k_1 + k_2 + l_1 + l_2 = 2r$, $1 \le l_i \le r$,
$$
\left<\tau_{k_1}h^{l_1}, \tau_{k_2}h^{l_2}\right>_d =
\frac{(-1)^{d(r + 1) + l_1 + k_2 + 1}}{d^{k_1 + k_2 + 1}} C^{2r -
(l_1 + l_2)}_{r - l_1},
$$
and other descendent invariants vanish. In particular, the only
non-trivial two-point function without descendents in degree
$d\ell$ is given by
$$
\left<h^r, h^r\right>_d = (-1)^{(d - 1)(r + 1)}\frac{1}{d}.
$$
\end{proposition}

\Proof We consider the invariant without descendents first. Since
the virtual dimension is $2r$, only $\langle h^r, h^r \rangle_d$
survives. Using the above equivalence relations, we may decrease
the power of $e_1^*h$ one by one. In each step only the second
term in the resulting three terms has nontrivial contribution.
Indeed, for the first term any addition to the power of $e_2^*h^r$
leads to zero.

For the third boundary splitting terms, write $[\Delta(X)] =
\sum_i T^i \otimes T_i$. For each $i$, since $\dim X = 2r + 1$ one
of $T^i$ or $T_i$ must have degree strictly bigger than $r$. If
$\beta_1 = d_1\ell$, $\beta_2 = d_2\ell$ with $d_i \ne 0$ then one
of the integral, hence the product, must vanish by the vanishing
lemma.

This is what happens now. We apply the divisor relation to $i = 1$
and $j = 2$. Since $n = 2$, we find $n_1 = |A| = 1$, $n_2 = |B| =
1$ and in the splitting we have sum of product of two-point
invariants. The degree in each side is non-zero since there is no
constant genus zero stable map with two marked points. So the
splitting terms vanish.

We apply the divisor relation repeatedly to compute
$$
\langle h^r, h^r\rangle_d = d\langle h^{r - 1}, \tau_1
h^r\rangle_d = \cdots = d^{r - 1}\langle h, \tau_{r - 1}
h^r\rangle_d = d^r\langle \tau_{r - 1} h^r\rangle
$$
where the last equality is by the divisor axiom. Now we plug in
Corollary~\ref{one-point} with $(k, l) = (r - 1, r)$ and the
statement follows.

For descendent invariants we proceed in the same manner. For
simplicity we abuse the notation by denoting $\left<\cdots, \psi^s
\alpha, \cdots\right>_\beta = \left<\cdots, \tau_s \alpha,
\cdots\right>_\beta$. Let $s \ge 1$, $l + m + s = 2r$ and consider
\begin{align*}
\left<h^l, \psi^s h^m\right>_d &= \left<h^{l - 1}, \psi^s h^{m +
1}\right>_d
 + (h, d\ell)\left<h^{l - 1}, \psi^{s + 1} h^m\right>_d\\
&= \left<h^{l - 1}, (h + d\psi)\psi^s h^m\right>_d = \cdots\\
&= \left<h, (h + d\psi)^{l - 1}\psi^s h^m\right>_d.
\end{align*}
Notice that the splitting terms are all zero as before. Now the
divisor axiom of descendent invariants gives
$$
d\left<(h + d\psi)^{l - 1}\psi^s h^m\right>_d + \left<(h +
d\psi)^{l - 1}\psi^{s - 1} h^{m + 1}\right>_d,
$$
which leads to the reduction formula:
$$
\left<h^l, \psi^s h^m\right>_d = \left<(h + d\psi)^l\psi^{s - 1}
h^m\right>_d.
$$

Notice that this equals the constant term in $z$ in
\begin{align*}
&\left<\sum_{k \ge 0}\frac{\psi^k}{z^k}(h + dz)^l z^{s - 1}
h^m\right>_d \\
&= z^{s + 1}e_{1*}\left(\frac{e(U_d)}{z(z - \psi)}.e_1^*\left((h +
dz)^l.h^m\right)\right)\\
&= (-1)^{(d - 1)(r + 1)}z^{s + 1}(h + dz)^{l - (r + 1)}.h^m \\
&= (-1)^{(d - 1)(r + 1)}\frac{z^{r - m}}{d^{(r + 1) - l}} \left(1
+ \frac{h}{dz}\right)^{l - (r + 1)}.h^m,
\end{align*}
which is
$$
\frac{(-1)^{d(r + 1) + r + 1}}{d^{r + 1 - l + r - m}}\, C^{l - (r
+ 1)}_{r - m} = \frac{(-1)^{d(r + 1) + l + s + 1}}{d^{s + 1}}\,
C^{2r - (l + m)}_{r - m}.
$$

In general from $\psi_1 = -\psi_2 + [D_{1|2}]^{virt}$, we find
$$
\left<\tau_{k_1}h^{l_1}, \tau_{k_2}h^{l_2}\right>_d =
-\left<\tau_{k_1 - 1}h^{l_1}, \tau_{k_2 + 1}h^{l_2}\right>_d =
\cdots = (-1)^{k_1} \left<h^{l_1}, \tau_{k_1 + k_2}
h^{l_2}\right>_d
$$
since the splitting terms all vanishes. The result follows.
\Endproof

For $n \ge 3$, it is known that for any three different markings
$i$, $j$ and $k$, $\psi_j = [D_{ik|j}]^{virt}$. By plugging this
into (\ref{e:rec1}), we get
$$
e_i^*L = e_j^*L + \sum_{\beta_1 + \beta_2 = \beta}((\beta_2.L)
[D_{ik, \beta_1|j, \beta_2}]^{virt} - (\beta_1.L) [D_{i,
\beta_1|jk, \beta_2}]^{virt}).
$$
In our special case this reads as
$$
e_i^*h = e_j^*h + \sum_{d_1 + d_2 = d}(d_2 [D_{ik, d_1|j,
d_2}]^{virt} - d_1 [D_{i, d_1|jk, d_2}]^{virt}).
$$
Notice that now $d_i$ is allowed to be zero.

\begin{lemma}
For $n \ge 3$,
\begin{align*}
&\langle h^{l_1 + 1}, h^{l_2}, h^{l_3}, \ldots \rangle_{n, d}\\
&= \langle h^{l_1}, h^{l_2 + 1}, h^{l_3},\ldots\rangle_{n, d} +
d\langle h^{l_1 + l_3}, h^{l_2},\ldots\rangle_{n -1, d} -d \langle
h^{l_1}, h^{l_2 + l_3}, \ldots\rangle_{n - 1, d}.
\end{align*}
Note that for $l_1 = 0$ this recovers the divisor axiom.
\end{lemma}

\Proof As in the previous theorem, the boundary terms with
non-trivial degree must vanish. For degree zero, the only
non-trivial invariants are three-point functions, hence we are
left with
\begin{align*}
&\langle h^{l_1 + 1}, h^{l_2}, h^{l_3}, \ldots \rangle_{n, d}\\
&= \langle h^{l_1}, h^{l_2 + 1}, h^{l_3},\ldots\rangle_{n, d}\\
&\quad + \sum_i d\langle h^{l_1}, h^{l_3}, T_i\rangle_0 \langle
T^i, h^{l_2},\ldots\rangle_d  - \sum_i d\langle T^i,
h^{l_1},\ldots\rangle_d \langle h^{l_2}, h^{l_3}, T_i\rangle_0.
\end{align*}
For the first boundary sum, in the diagonal decomposition
$[\Delta(X)] = \sum T_i\otimes T^i$ we may choose basis so that
$h^{l_1 + l_3}$ appear in $\{T^i\}$. Then the above degree zero
invariants survive only in one term which is equal to 1. The same
argument applies to the second sum too. So the above expression
equals
$$
\langle h^{l_1}, h^{l_2 + 1}, h^{l_3},\ldots\rangle_{d} + d\langle
h^{l_1 + l_3}, h^{l_2},\ldots\rangle_{d} -d \langle h^{l_1},
h^{l_2 + l_3}, \ldots\rangle_{d}
$$
as expected. \Endproof

In light of \eqref{e:virt} and results above, Theorem~\ref{main}
can be reformulated as the following equation
\begin{equation} \label{e:main}
\langle h^{l_1}, h^{l_2},\ldots, h^{l_n}\rangle_d = (-1)^{(d -
1)(r + 1)}N_{l_1,\ldots,l_n}d^{n -3},
\end{equation}
which we will now prove.

\Proof (of \eqref{e:main}, or equivalently Theorem~\ref{main}.)

We will prove the theorem by induction on $n \in \mathbb{N}$. The
case $n \le 2$ are already proven before. We treat the case $n =
3$ first.

Consider $\langle h^{l_1}, h^{l_2}, h^{l_3} \rangle_d$ with $l_1 +
l_2 + l_3 = 2r + 1$ and $l_1 \le l_2 \le l_3$. If $l_1 = 1$ then
$l_2 = l_3 = r$ and so
$$
\langle h, h^r, h^r \rangle_d  = d\langle h^r, h^r \rangle_d =
(-1)^{(d - 1)(r + 1)}.
$$
If $l_1 \ge 2$, then $l_2 \le r - 1$ and
$$
\langle h^{l_1}, h^{l_2}, h^{l_3} \rangle_d = \langle h^{l_1- 1},
h^{l_2 + 1}, h^{l_3} \rangle_d + d\langle h^{l_1 + l_3 - 1},
h^{l_2}\rangle_d - d\langle h^{l_1 - 1}, h^{l_2 + l_3}\rangle_d.
$$
But then both $l_1 + l_3 - 1$ and $l_2 + l_3$ are larger than $r +
1$ and the boundary terms vanish individually. By reordering
$l_2$, $l_3$ if necessary, and repeating this procedure we are
reduced to the case $l_1 = 1$ and proof for $n = 3$ is completed.

Suppose the theorem holds up to $n - 1$ (with $n \ge 4$). The
above lemma and the induction hypothesis imply that
\begin{align*}
&\langle h^{l_1}, h^{l_2}, h^{l_3}, \ldots \rangle_d\\
&= \langle h^{l_1 - 1}, h^{l_2 + 1}, h^{l_3},\ldots\rangle_d +
d\langle h^{l_1 + l_3 - 1}, h^{l_2},\ldots\rangle_d - d \langle
h^{l_1 - 1}, h^{l_2 + l_3}, \ldots\rangle_d\\
&= \langle h^{l_1 - 1}, h^{l_2 + 1}, h^{l_3},\ldots\rangle_d +
(N_{l_1 + l_3 - 1,l_2,\ldots} - N_{l_1 - 1,l_2 + l_3,\ldots})d^{n
- 3}.
\end{align*}
By repeating this procedure, $l_1$ is decreased to one and we get
$$
\langle h^{l_1}, h^{l_2},\ldots, h^{l_n}\rangle_d = (-1)^{(d -
1)(r + 1)}N_{l_1,\ldots,l_n}d^{n -3},
$$
where $N_{l_1,\ldots,l_n}$ is given by $N_*$'s in one lower level.
The proof is complete. \Endproof

Similar methods apply to descendent invariants:

\begin{proposition}\label{main-d}
The only three-point descendent invariants of extremal classes
$d\ell$, up to permutations of insertions, are given by
$$
\left<h^{l_1}, h^{l_2}, \tau_{k_3} h^{l_3}\right>_d =
\frac{(-1)^{d(r + 1) + l_3 + 1}}{d^{k_3}}\,C^{k_3 + 1}_{r - (l_1 +
l_2)},
$$
where $l_1 + l_2 + l_3 + k_3 = 2r + 1$ and by convention $C^m_n =
0$ if $n < 0$.

More generally, an $n$-point descendent invariant $\langle
\prod_{i = 1}^n \tau_{k_i} h^{l_i} \rangle_d$ with $n \ge 3$ is
non-zero only if there are at least two insertions being free of
descendents, say $k_1 = k_2 = 0$. In such cases, there are
universal constants $N_{k, l} \in \mathbb{Z}$ such that
$$
\left<h^{l_1}, h^{l_2}, \tau_{k_3} h^{l_3},\ldots, \tau_{k_n}
h^{l_n}\right>_d = N_{k, l}\,d^{n - 3 - \sum k_i}.
$$
\end{proposition}

\Proof Let $n \ge 3$ and assume that $0 \le k_1 \le k_2 \le \cdots
\le k_n$. If $k_2 \ge 1$ then use $\psi_2 = [D_{2|13}]^{virt}$ we
get
$$
\left<\tau_{k_1} h^{l_1}, \ldots, \tau_{k_n} h^{l_n} \right>_d =
\sum_{i;\, d_1 + d_2 = d} \left<\tau_{k_2 - 1} h^{l_2},\cdots,
T_i\right>_{d_1} \left<T^i, \tau_{k_1} h^{l_1}, \tau_{k_3}
h^{l_3}, \cdots\right>_{d_2}.
$$
We separate two cases. If the first factor is a two-point function
then it is non-zero only if $T_i = h^j$ for some $j \le r$. But
then $\deg T^i > r$ and the right factor vanishes since it
contains $\psi$ classes. For other cases, both factors contain
$\psi$ classes hence the factor with $\deg T_i > r$ (or $\deg T^i
> r$) must vanish.

For three-point invariants, from $\psi_3 = [D_{3|12}]^{virt}$ we
get as before that
\begin{align*}
\left<h^{l_1}, h^{l_2}, \tau_{k_3} h^{l_3}\right>_d &= \sum_{i;\,
d_1 + d_2 = d} \left<\tau_{k_3 - 1} h^{l_3}, T_i\right>_{d_1}
\left<T^i, h^{l_1}, h^{l_2}\right>_{d_2}\\
&= \left<\tau_{k_3 - 1} h^{l_3}, h^{l_1 + l_2}\right>_d
\end{align*}
and the formula follows from the two-point case.

Similarly, for $n \ge 4$, if $k_i \ne 0$ then from $\psi_i =
[D_{i|12}]^{virt}$ we get
$$
\left<h^{l_1}, h^{l_2}, \ldots, \tau_{k_i} h^{l_i},
\ldots\right>_d = \left<\tau_{k_i - 1} h^{l_i}, h^{l_1 + l_2},
\ldots\right>_d.
$$
The result follows from an induction on $n$. \Endproof

\section{Degeneration analysis}  \label{s:4}

Our next task is to \emph{compare} the genus zero Gromov--Witten
invariants of $X$ and $X'$ for curve classes other than the
flopped curve. Naively, one may wish to ``decompose'' the
varieties into the neighborhoods of exceptional loci and their
complements. As the latter's are obviously isomorphic, one is
reduced to study the local case. The degeneration formula
\cite{LiRu} \cite{Li} \cite{IoPa} provides a rigorous formulation
of the above naive picture.

\subsection{The degeneration formula}
Our presentation of degeneration formula below mostly follows that
of \cite{Li} and \cite{LiuYau}. We have, however, chosen to use
the ``numerical form'' rather than the ``cycle form'' in the
exposition.

Given a relative pair $(Y, E)$ with $E \hookrightarrow Y$ a smooth
divisor, the relative Gromov--Witten invariants are defined in the
following way. Let $\Gamma = (g, n, \beta, \rho, \mu)$ with $\mu =
(\mu_1, \ldots, \mu_\rho) \in \mathbb{N}^\rho$ a partition of the
intersection number $(\beta.E)= |\mu| := \sum_{i = 1}^\rho \mu_i$.
For $A \in H^*(Y)^{\otimes n}$ and $\s \in H^*(E)^{\otimes \rho}$,
the relative invariant of stable maps with topological type
$\Gamma$ (i.e.~with contact order $\mu_i$ in $E$ at the $i$-th
contact point) is
$$
\langle A \mid \s, \mu \rangle^{(Y, E)}_{\Gamma} :=
\int_{[\overline{M}_{\Gamma}(Y, E)]^{virt}} e_Y^* A \cup e_E^* \s
$$
where $e_Y: \overline{M}_{\Gamma}(Y, E) \to Y^n$, $e_E:
\overline{M}_{\Gamma}(Y, E) \to E^\rho$ are evaluation maps on
marked points and contact points respectively.

If $\Gamma = \coprod_\pi \Gamma^\pi$, the relative invariants
(with disconnected domain curves)
$$
\langle A \mid \s, \mu \rangle^{\bullet (Y, E)}_{\Gamma} :=
\prod\nolimits_\pi \langle A \mid \s, \mu \rangle^{(Y,
E)}_{\Gamma^\pi}
$$
are defined to be the product of the connected components.

We apply the degeneration formula to the following situation.
Let $X$ be a smooth variety and $Z \subset X$ be a smooth subvariety.
Let $\Phi:W \to \mathscr{X}$ be its \emph{degeneration to the
normal cone}, the blow-up of $X \times \mathbb{A}^1$ along $Z
\times \{0\}$. Denote $t\in\mathbb{A}^1$ the deformation
parameter. Then $W_t \cong X$ for all $t\ne 0$ and $W_0 = Y_1 \cup
Y_2$ with
$$
   \phi =\Phi|_{Y_1}: Y_1 \to X
$$
the blow-up along $Z$ and
$$
  p =\Phi|_{Y_2}: Y_2 := \mathbb{P}_Z(N_{Z/X}\oplus \mathscr{O})
  \to Z\subset X
$$
the projective completion of the normal bundle. $Y_1 \cap Y_2 =: E
= \mathbb{P}_Z(N_{Z/X})$ is the $\phi$-exceptional divisor which
consists of ``the infinity part'' of the projective bundle
$\mathbb{P}_Z(N_{Z/X}\oplus \mathscr{O})$.

\begin{figure}
\begin{center}
\includegraphics[width=3.5in]{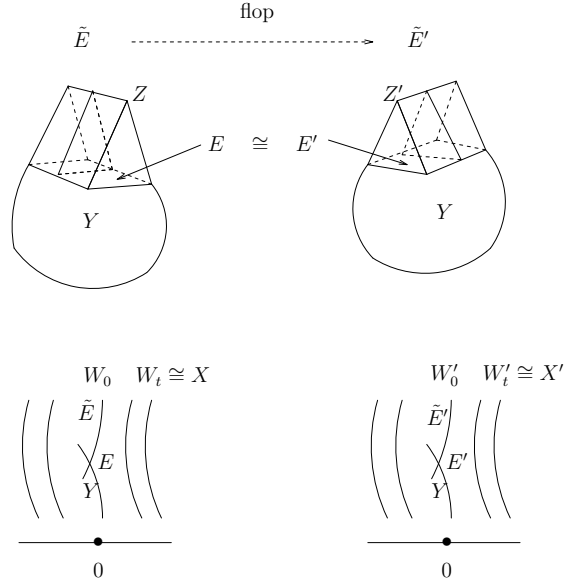}
\caption{Degeneration to the normal cone for ordinary flops.}
\end{center}
\end{figure}

Since the family $W \to \mathbb{A}^1$ is a degeneration of a
trivial family, all cohomology classes $\alpha\in
H^*(X,\mathbb{Z})^{\oplus n}$ have global liftings and the
restriction $\alpha(t)$ on $W_t$ is defined for all $t$. Let $j_i:
Y_i\hookrightarrow W_0$ be the inclusion maps for $i = 1, 2$. Let
$\{\e_i\}$ be a basis of $H^*(E)$ with $\{\e^i\}$ its dual basis.
$\{\e_I\}$ forms a basis of $H^*(E^\rho)$ with dual basis
$\{\e^I\}$ where $|I| = \rho$, $\e_I = \e_{i_1}\otimes \cdots
\otimes \e_{i_\rho}$. The \emph{degeneration formula} expresses
the absolute invariants of $X$ in terms of the relative invariants
of the two smooth pairs $(Y_1,E)$ and $(Y_2, E)$:
$$
\langle\alpha\rangle_{g,n,\beta}^X = \sum_{I} \sum_{\eta\in
\Omega_\beta} C_\eta \left.\Big\langle j_1^*\alpha(0) \,\right|\,
\e_I, \mu \Big\rangle_{\Gamma_1}^{\bullet (Y_1,E)}
\left.\Big\langle j_2^*\alpha(0) \,\right|\, \e^I, \mu
\Big\rangle_{\Gamma_2}^{\bullet (Y_2,E)}.
$$
Here $\eta = (\Gamma_1, \Gamma_2, I_\rho)$ is an {\it admissible
triple} which consists of (possibly disconnected) topological
types
$$
\Gamma_i = \coprod\nolimits_{\pi = 1}^{|\Gamma_i|} \Gamma_i^\pi
$$
with the same partition $\mu$ of contact order under the
identification $I_\rho$ of contact points. The gluing $\Gamma_1
+_{I_\rho} \Gamma_2$ has type $(g, n, \beta)$ and is connected. In
particular, $\rho = 0$ if and only if that one of the $\Gamma_i$
is empty. The total genus $g_i$, total number of marked points
$n_i$ and the total degree $\beta_i \in NE(Y_i)$ satisfy the
splitting relations
\[
 \begin{split}
 g &= g_1 + g_2 + \rho + 1 - |\Gamma_1| -|\Gamma_2|, \\
 n &= n_1 + n_2, \\
 \beta &= \phi_*\beta_1 + p_*\beta_2.
\end{split}
\]

The constants $C_\eta = m(\mu)/|{\rm Aut}\,\eta|$, where $m(\mu) =
\prod \mu_i$ and ${\rm Aut}\,\eta = \{\, \sigma \in S_\rho \mid
\eta^\sigma = \eta \,\}$. (When a map is decomposed into two
parts, an (extra) ordering to the contact points is assigned. The
automorphism of the decomposed curves will also introduce an extra
factor. These contribute to ${\rm Aut}\,\eta$.) We denote by
$\Omega$ the set of equivalence classes of all admissible triples;
by $\Omega_\beta$ and $\Omega_\mu$ the subset with fixed degree
$\beta$ and fixed contact order $\mu$ respectively.

Given an ordinary flop $f : X \dashrightarrow X'$, we apply
degeneration to the normal cone to both $X$ and $X'$. Then $Y_1
\cong Y'_1$ and $E = E'$, by the definition of ordinary flops. The
following notations will be used
\[
 Y := {\rm Bl}_Z X \cong Y_1 \cong Y'_1, \quad
 \tilde{E} :=  \mathbb{P}_Z(N_{Z/X}\oplus \mathscr{O}), \quad
 \tilde{E}' :=  \mathbb{P}_{Z'}(N_{Z'/X'}\oplus \mathscr{O}).
\]

\begin{remark}
For simple $\mathbb{P}^r$ flops, $Y_2 \cong
\mathbb{P}_{\mathbb{P}^r}(\mathscr{O}(-1)^{\oplus(r + 1)}\oplus
\mathscr{O}) \cong Y'_2$. However the gluing maps of $Y_1$ and
$Y_2$ along $E$ for $X$ and $X'$ differ by a twist which
interchanges the order of factors in $E = \mathbb{P}^r \times
\mathbb{P}^r$. Thus $W_0 \not\cong W_0'$ and it is necessary to
study the details of the degenerations. In general, $f$ induces an
ordinary flop $\tilde f: Y_2 \dashrightarrow Y_2'$ of the same
type which is the local model of $f$.
\end{remark}

\subsection{Liftings of cohomology insertions}

Next we discuss the presentation of $\alpha(0)$. Denote by
$\iota_1 \equiv j: E \hookrightarrow Y_1 = Y$ and $\iota_2: E
\hookrightarrow Y_2 = \tilde E$ the natural inclusions. The class
$\alpha(0)$ can be represented by $(j_1^*\alpha(0),
j_2^*\alpha(0)) = (\alpha_1, \alpha_2)$ with $\alpha_i \in
H^*(Y_i)$ such that
\begin{equation} \label{e:4.2.1}
\iota_1^*\alpha_1 = \iota_2^*\alpha_2 \quad\mbox{and}\quad
\phi_*\alpha_1 + p_*\alpha_2 = \alpha.
\end{equation}
Such representatives are called {\it liftings} which are by no
means unique. The flexibility on different choices will be useful.

One choice of the lifting is
\begin{equation} \label{e:4.2.2}
\alpha_1 = \phi^*\alpha \quad \text{and} \quad \alpha_2 =
p^*(\alpha|_{Z}),
\end{equation}
since they satisfy the conditions \eqref{e:4.2.1}:
$(\alpha_1,\alpha_2)$ restrict to the same class in $E$
and push forward to $\alpha$ and $0$ in $X$ respectively.
More generally:

\begin{lemma}\label{switch}
Let $\alpha(0) = (\alpha_1, \alpha_2)$ be a choice of lifting.
Then
$$
\alpha(0) = (\alpha_1 - \iota_{1*} e, \alpha_2 + \iota_{2*} e)
$$
is also a lifting for any class $e$ in $E$ of the same dimension
as $\alpha$. Moreover, any two liftings are related in this
manner. In particular, $\alpha_1$ and $\alpha_2$ are uniquely
determined by each other.
\end{lemma}

\Proof The first statement follows from the facts that
$$
\iota_1^*\iota_{1*} e = (e.c_1(N_{E/Y}))_E = -(e.c_1(N_{E/\tilde
E}))_E = -\iota_2^*\iota_{2*} e
$$
and $- \phi_*\iota_{1*}e + p_*\iota_{2*}e = 0$ (since $\phi \circ
\iota_1 = p \circ \iota_2 = \bar \phi: E \to Z$).

For the second statement, let $(\alpha_1, \alpha_2)$ and $(a_1,
a_2)$ be two liftings. From
$$
\phi_*(\alpha_1 - a_1) = -p_*(\alpha_2 - a_2) \in H^*(Z),
$$
we have that $\phi^*\phi_*(\alpha_1 - a_1)$ is a class in $E$.
Hence $\alpha_1 - a_1 = \iota_{1*}e$ for $e \in H^*(E)$. It
remains to show that if $(a_1, a_2)$ and $(a_1, \tilde a_2)$ are
two liftings then $a_2 = \tilde a_2$. Indeed by \eqref{e:4.2.1},
$\iota_2^*(a_2 - \tilde a_2) = 0$. Hence by Lemma \ref{l:4-3}
below $z := a_2 - \tilde a_2 \in i_* H^*(Z)$. By \eqref{e:4.2.1}
again $z = p_*z = p_*(a_2 - \tilde a_2) = 0$. \Endproof

For an ordinary flop $f:X\dashrightarrow X'$, we compare the
degeneration expressions of $X$ and $X'$.
For a given admissible triple $\eta = (\Gamma_1, \Gamma_2, I_\rho)$ on
the degeneration of $X$, one may pick the corresponding
$\eta' = (\Gamma_1', \Gamma_2', I'_\rho)$ on the degeneration of $X'$
such that $\Gamma_1 = \Gamma_1'$.
Since
$$
\phi^*\alpha - \phi'^*\T \alpha \in \iota_{1*} H^*(E) \subset H^*(Y),
$$
Lemma~\ref{switch} implies that one can choose
$\alpha_1=\alpha'_1$. This can be done, for example, by modifying
the choice of \eqref{e:4.2.2} $j_1^*\alpha(0) = \phi^*\alpha$ and
$j_1'^*\T\alpha(0) = \phi'^* \T\alpha$ by adding suitable classes
in $E$ to make them equal. The above procedures identify relative
invariants on the $Y_1 = Y = Y'_1$ from both sides term by term,
and we are left with the comparison of the corresponding relative
invariants on $\tilde E$ and $\tilde E'$. The following simple
lemma is useful.

\begin{lemma} \label{l:4-3}
Let $\tilde E = \mathbb{P}_Z(N \oplus \mathscr{O})$ be a
projective bundle with base $i: Z \hookrightarrow \tilde E$ and
infinity divisor $\iota_2: E = \mathbb{P}_Z(N) \hookrightarrow
\tilde E$. Then the kernel of the restriction map $\iota_2^* :
H^*(\tilde E) \to H^*(E)$ is $i_* H^*(Z)$.
\end{lemma}

\Proof
$i_* H^*(Z)$ obviously lies in the kernel of $\iota_2^*$.
The fact it is the entire kernel can be seen, for example,
by a dimension count.
\Endproof

The ordinary flop $f$ induces an ordinary flop
\[
\tilde f : \tilde E \dashrightarrow \tilde E'
\]
on the local model. Moreover $\tilde{f}$ may be considered as a
family of simple ordinary flops $\tilde f_t: \tilde E_t
\dashrightarrow \tilde E'_t$ over the base $S$, where $t \in S$
and $\tilde E_t$ is the fiber of $\tilde E \to Z \to S$ etc..
Denote again by $\T$ the cohomology correspondence induced by the
graph closure. Then

\begin{proposition}[Cohomology reduction to local models] \label{coh-red}
Let $f:X\dashrightarrow X'$ be a $\mathbb{P}^r$ flop over base $S$
with $\dim S = s$. Let $\alpha \in H^*(X)$ with liftings
$\alpha(0) = (\alpha_1, \alpha_2)$ and $\T\alpha(0) = (\alpha'_1,
\alpha'_2)$. Then
$$
\alpha_1 = \alpha'_1 \quad \Longleftrightarrow \quad \T\alpha_2 =
\alpha'_2.
$$
\end{proposition}

\Proof Let $\alpha \in H^{2l}(X)$ with $l \in \frac{1}{2}
\mathbb{N}$. If $l
> \dim Z = r + s$ then $\alpha|_Z =0$. By \eqref{e:4.2.2} and
Lemma \ref{switch}, all liftings take the form $\alpha(0) =
(\alpha - \iota_{1*}e, \iota_{2*}e)$ and $\T\alpha(0) = (\alpha -
\iota'_{1*}e', \iota'_{2*}e')$ for $e, e'$ being classes in $E$.
In this case the proof is trivial since $\T$ is the identity map
on $H^*(E)$. So we may assume that $l \le r + s$.

($\Rightarrow$) From the contact order condition
$\iota_2^*\alpha_2 = \iota_1^*\alpha_1 = \iota_1'^*\alpha'_1 =
\iota_2'^*\alpha'_2$ and the fact that $\tilde f$ is an
isomorphism outside $Z$, we get
$$
\iota_2'^*(\T \alpha_2 - \alpha'_2) = \T \iota_2^*\alpha_2 -
\iota_2'^* \alpha'_2 = \iota_2^*\alpha_2 - \iota_2'^* \alpha'_2 =
0.
$$
Thus $\T \alpha_2 - \alpha'_2 = i'_*z'$ for some $z' \in H^{2(l -
(r + 1))}(Z')$ (where $i': Z' \hookrightarrow \tilde{E}'$) by
Lemma~\ref{l:4-3} and the fact that ${\rm codim}_{\tilde E'} Z' =
r + 1$.

For simple flops, $s = 0$ and then $l - (r + 1) \le s - 1 < 0$. So
$z' = 0$ and we are done. In general we restrict the equation to
each fiber $\tilde f_t: \tilde E_t \to \tilde E'_t$. Since
$\bar\Gamma_{\tilde f}|_t = \bar\Gamma_{\tilde f_t}$, by the case
of simple flops we get $(\T \alpha_2 - \alpha'_2)|_{\tilde E'_t} =
0$ for all $t \in S$. That is, $z'$ is a class supported in the
fiber of $p': Z' \to S$. But then ${\rm codim}_{\tilde E'}z' \ge s
+ r + 1
> l$, which implies that $z' =0$.

($\Leftarrow$) For ease of notations we omit the imbedding maps of
$E$ into $Y$, $\tilde E$ and $\tilde E'$. By \eqref{e:4.2.2} and
Lemma \ref{switch} we have $\alpha_1 = \phi^*\alpha - e_1$ and
$\alpha_1' = \phi'^*\T\alpha - e_1'$ for some classes $e_1$,
$e_1'$ in $E$. Thus $\alpha_1' = \alpha_1 - e$ for some class $e$
in $E$. By Lemma \ref{switch} again $\alpha(0)$ has a lifting
$(\alpha_1 - e, \alpha_2 + e) = (\alpha_1', \alpha_2 + e)$ and by
the first part of this proposition we must have $\T(\alpha_2 + e)
= \alpha_2'$. By assumption $\T\alpha_2 = \alpha_2'$, hence $\T e
= 0$ and then $e = 0$. \Endproof

\begin{remark} \label{r:4.5}
Proposition~\ref{coh-red} (with cohomology groups being replaced
by Chow groups) leads to an alternative proof of equivalence of
Chow motives under ordinary flops. Indeed the equivalence of Chow
groups for simple flops is easy to establish. The degeneration to
normal cone then allows us to reduce the general case to the local
case and then to the local simple case.
\end{remark}

\subsection{Reduction to relative local models}
First notice that $A_1(\tilde E) = \iota_{2*} A_1(E)$ since both
are projective bundles over $Z$. We then have
$$
\phi^*\beta = \beta_1 + \beta_2
$$
by regarding $\beta_2$ as a class in $E \subset Y$. Indeed
$\phi_*(\beta_1 + \beta_2) = \phi_*\beta_1 + p_*\beta_2 = \beta$
and
$$
((\beta_1 + \beta_2).E)_Y = (\beta_1.E)_Y - (\beta_2.E)_{\tilde E}
= |\mu| - |\mu| = 0
$$
(where $N_{E/\tilde E} \cong N_{E/Y}^*$ is used). These
characterize the class $\phi^*\beta$.

We consider only the case $g = 0$.
Define the generating series
$$
\langle A \mid \s, \mu \rangle^{(\tilde{E}, E)}
 := \sum_{\beta_2 \in NE(\tilde E)}
 \frac{1}{|{\rm Aut}\,\mu|} \langle A \mid \s, \mu
 \rangle_{\beta_2}^{(\tilde E, E)}\,\q^{\beta_2}.
$$
and the similar one with possibly disconnected domain curves
$$
\langle A \mid \s, \mu \rangle^{\bullet (\tilde E, E)} :=
\sum_{\Gamma;\, \mu_\Gamma = \mu} \frac{1}{|{\rm Aut}\,\Gamma|}
\langle A \mid \s, \mu \rangle_{\Gamma}^{\bullet (\tilde E,
E)}\,\q^{\beta^\Gamma}.
$$

\begin{proposition}\label{reduction}
To prove $\T \langle \alpha \rangle^X \cong \langle \T\alpha
\rangle^{X'}$ (for all $\alpha$), it is enough to show that
\begin{equation} \label{e:4.3.1}
\T\langle A \mid \s, \mu \rangle^{(\tilde{E}, E)} \cong
\langle \T A \mid \s, \mu \rangle^{(\tilde{E}', E)}
\end{equation}
for all $A, \s, \mu$.
\end{proposition}

\Proof For the $n$-point function $\langle \alpha \rangle^X =
\sum_{\beta \in NE(X)} \langle \alpha \rangle^X_\beta\, \q^\beta$,
the degeneration formula gives
\begin{align*}
\langle \alpha \rangle^X &= \sum_{\beta \in NE(X)}\sum_{\eta \in
\Omega_\beta} \sum_I C_\eta \langle \alpha_1 \mid \e_I, \mu
\rangle_{\Gamma_1}^{\bullet (Y_1,E)} \langle \alpha_2 \mid
\e^I, \mu \rangle_{\Gamma_2}^{\bullet (Y_2,E)}\,\q^{\phi^*\beta} \\
&= \sum_{\mu} \sum_I \sum_{\eta \in \Omega_\mu} C_\eta
\left(\langle \alpha_1 \mid \e_I, \mu \rangle_{\Gamma_1}^{\bullet
(Y_1,E)}\,\q^{\beta_1}\right) \left(\langle \alpha_2 \mid \e^I,
\mu \rangle_{\Gamma_2}^{\bullet (Y_2,E)}\,\q^{\beta_2}\right).
\end{align*}

To simplify the generating series, we consider also absolute
invariants $\langle \alpha \rangle^{\bullet X}$ with possibly
disconnected domain curves as before. Then by comparing the order
of automorphisms,
$$
\langle \alpha \rangle^{\bullet X} = \sum_{\mu} m(\mu) \sum_I
\langle \alpha_1 \mid \e_I, \mu \rangle^{\bullet (Y_1,E)} \langle
\alpha_2 \mid \e^I, \mu \rangle^{\bullet (Y_2,E)}.
$$

To compare $\T \langle \alpha \rangle^{\bullet X}$ and $\langle \T
\alpha \rangle^{\bullet X'}$, by Proposition~\ref{coh-red} we may
assume that $\alpha_1 = \alpha_1'$ and $\alpha_2' = \T\alpha_2$.
This choice of cohomology liftings identifies the relative
invariants of $(Y_1, E)$ and those of $(Y'_1, E')$ with the same
topological types. It remains to compare
$$
\langle \alpha_2 \mid \e^I, \mu \rangle^{\bullet (\tilde E, E)}
\quad \mbox{and} \quad \langle \T\alpha_2 \mid \e^I, \mu
\rangle^{\bullet (\tilde E', E)}.
$$

We further split the sum into connected invariants. Let
$\Gamma^\n$ be a connected part with the contact order
$\mu^\n$ induced from $\mu$.
Denote $P: \mu = \sum_{\n \in P} \mu^\n$ a partition of $\mu$
and $P(\mu)$ the set of all such partitions. Then
$$
\langle A \mid \s, \mu \rangle^{\bullet (\tilde E, E)} = \sum_{P
\in P(\mu)} \prod_{\n \in P} \sum_{\Gamma^\n} \frac{1}{|{\rm
Aut}\,\mu^\n|} \langle A^\n \mid \s^\n, \mu^\n
\rangle_{\Gamma^\n}^{(\tilde E, E)}\,\q^{\beta^{\Gamma^\n}}.
$$

If one fixes the above data in the summation of \eqref{e:4.3.1},
then the only index to be summed over is $\beta^{\Gamma^\n}$ on
$\tilde{E}$. This reduces the problem to $\langle A^\n \mid \s^\n,
\mu^\n \rangle^{(\tilde{E}, E)}$. \Endproof

\begin{remark}
Here is a brief comment on the term
$$
\T \langle \alpha_2 \mid \e^I, \mu \rangle^{(\tilde E, E)} =
\sum_{\beta_2 \in NE(\tilde E)} \frac{1}{|{\rm Aut}\,\mu|} \langle
\alpha_2 \mid \e^I, \mu \rangle_{\beta_2}^{(\tilde E, E)}\,\q^{\T
\beta_2}.
$$
Since $\tilde E$ is a projective bundle, $NE(\tilde E) = i_* NE(Z)
\oplus \mathbb{Z}_+\gamma$ with $\gamma$ the fiber line class of
$\tilde E \to Z$. The point is that, for $\beta_2 \in NE(\tilde
E)$ it is in general not true that $\T \beta_2 \equiv \beta_2$ (in
$E$) is effective in $\tilde E'$.

Indeed, for simple ordinary flops, let $\gamma = \delta'$, $\delta
= \gamma'$ be the two line classes in $E \cong \mathbb{P}^r \times
\mathbb{P}^r$. It is easily checked that $\ell \sim \delta -
\gamma$ in $\tilde E$. Hence $\ell = -\ell'$ and $\gamma = \gamma'
+ \ell'$ and
$$
\beta_2 = d_1 \ell + d_2 \gamma = (d_2 - d_1)\ell' + d_2\gamma'.
$$
$\T\beta_2 \in NE(\tilde E')$ if and only if $d_2 \ge d_1$.
Therefore,
$$
\langle \alpha_2 \mid \e^I, \mu \rangle^{(\tilde E, E)} = \langle
\T \alpha_2 \mid \e^I, \mu \rangle^{(\tilde E', E)}
$$
cannot possibly hold term by term.
Analytic continuations are in general needed.
\end{remark}

\subsection{Relative to absolute}
Recall that we are now in the local relative case, with $X=\tilde{E}$.
We shall combine a method of Maulik and Pandharipande
(Lemma 4 in \cite{MP}) to further reduce the relative cases to the
absolute cases with at most descendent insertions along $E$.
Following \cite{MP}, we call the pair
$$
(\s, \mu) = \{(\s_1, \mu_1), \cdots, (\s_\rho, \mu_\rho)\}
$$
with $\s_i \in H^*(E)$, $\mu_i \in \mathbb{N}$ a \emph{weighted
partition}, a partition of contact orders weighted by cohomology
classes in $E$.

\begin{proposition} \label{p:4.8}
For an ordinary flop $\tilde E \dashrightarrow \tilde E'$, to prove
$$
\T\langle A \mid \s, \mu \rangle \cong \langle \T A \mid \s, \mu
\rangle
$$
for any $A$ and $(\s, \mu)$, it is enough to show that
$$
\T\langle A, \tau_{k_1} \s_1, \ldots, \tau_{k_\rho} \s_\rho
\rangle^{\tilde{E}} \cong \langle \T A, \tau_{k_1} \s_1, \ldots,
\tau_{k_\rho} \s_\rho \rangle^{\tilde{E}'}
$$
for any possible insertions $A \in H^*(\tilde E)^{\oplus n}$, $k_j
\in \mathbb{N}\cup \{0\}$ and $\s_j \in H^*(E)$. (Here we abuse
the notations and denote $\iota_{2*}\, \s \in H^*(\tilde E)$ by
the same symbol $\s$.)
\end{proposition}

The rest of this subsection is devoted to the proof of this
proposition which proceeds inductively on the triple ($|\mu|$,
$n$, $\rho$) in the lexicographical order with $\rho$ in the
reverse order. Given $\langle \alpha_1, \ldots, \alpha_n \mid \s,
\mu\rangle$, since $\rho \le |\mu|$, it is clear that there are
only finitely many triples of lower order. The proposition holds
for those cases by the induction hypothesis.

We apply degeneration to the normal cone for $Z \hookrightarrow
\tilde E$ to get $W \to \mathbb{A}^1$. Then $W_0 = Y_1 \cup Y_2$
with $\pi: Y_1 \cong \mathbb{P}_E(\mathscr{O}_E(-1, -1)\oplus
\mathscr{O}) \to E$ a $\mathbb{P}^1$ bundle and $Y_2 \cong \tilde
E$. Denote by $E_0 = E = Y_1 \cap Y_2$ and $E_\infty \cong E$ the
zero and infinity divisors of $Y_1$ respectively. The idea is to
analyze the degeneration formula for $\langle \alpha_1, \ldots,
\alpha_n, \tau_{\mu_1 - 1}\s_1, \ldots, \tau_{\mu_\rho - 1}
\s_\rho \rangle^{\tilde E}$. We follow the procedure used in the
proof of Proposition \ref{reduction} to split the generating
series of invariants with possibly disconnected domain curves,
according to the contact order. For $\beta = d_1 \ell + d_2 \gamma
\in NE(\tilde E)$, $c_1(\tilde E).\beta = d_2 c_1(\tilde
E).\gamma$, hence by the virtual dimension counting $d_2$ is
uniquely determined for a given generating series with fixed
cohomology insertions.

\begin{figure}
\begin{center}
\includegraphics[width=3.5in]{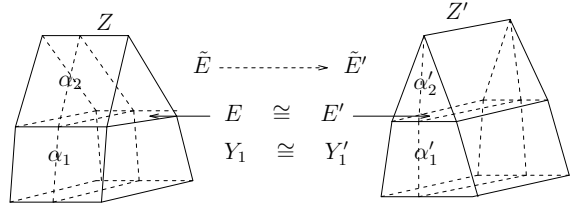}
\caption{Degeneration to normal cone for local models.}
\end{center}
\end{figure}

We observe that during the splitting of $\beta$'s, the ``main
terms'' with the highest total contact order only occur when the
curve classes in $Y_1$ are fiber classes. Indeed, let $(\beta_1,
\beta_2)$ be a splitting of $\beta$. Since
$$
NE(Y_1) = \mathbb{Z}_+\delta + \mathbb{Z}_+\bar\gamma +
\mathbb{Z}_+\gamma \quad \mbox{and}\quad  NE(Y_2) =
\mathbb{Z}_+\ell + \mathbb{Z}_+\gamma
$$
($\bar \gamma$ is the fiber class of $Y_1$), we have
$$
(\beta_1, \beta_2) = (a \delta + b \gamma + c \bar\gamma, d \ell +
e \gamma)
$$
subject to
$$
a, b, c, d, e \ge 0, \quad a + d = d_1,\quad c = d_2
$$
and the total contact order condition
$$
e = (\beta_2.E)_{\tilde E} = (\beta_1.E)_{Y_1} = -a - b + c.
$$
In particular, $e \le d_2$ with $e = d_2$ if and only if that $a =
b = 0$. In this case $\beta_1 = d_2\bar\gamma$ and the invariants
on $(Y_1, E)$ are fiber class integrals.

It is sufficient to consider $(\s_1, \ldots, \s_\rho) = \e_I =
(\e_{i_1}, \ldots, \e_{i_\rho})$. Since $\s_i|_Z = 0$, one may
choose the cohomology lifting $\s_i(0) = (\iota_{1*}\, \s_i, 0)$.
This ensures that insertions of the form $\tau_k\,\s$ must go to
the $Y_1$ side in the degeneration formula.

\begin{lemma}
For a general cohomology insertion $\alpha \in H^*(\tilde E)$, the
lifting can be chosen to be $\alpha(0) = (a, \alpha)$ for some
$a$.
\end{lemma}

\Proof $\alpha(0)$ may be chosen as $(\phi^* \alpha,
p^*(\alpha|_Z))$. Since $(\alpha - p^*(\alpha|_Z)).Z = 0$, the
class $e := \alpha - p^*(\alpha|_Z)$ can be taken to be supported
in $E$. Then Lemma \ref{switch} implies that $\alpha(0)$ can be
modified to be $(\phi^*\alpha - e, \alpha)$. \Endproof

From $\alpha(0) = (a, \alpha)$ and $\T\alpha(0) = (a', \T\alpha)$,
Lemma \ref{l:4-3} implies that $a = a'$. As before the relative
invariants on $(Y_1, E)$ can be regarded as constants under $\T$.
Then
\begin{align*}
&\langle \alpha_1, \ldots, \alpha_n, \tau_{\mu_1 - 1}\e_{i_1},
\ldots, \tau_{\mu_\rho - 1} \e_{i_\rho} \rangle^{\bullet\tilde E}
= \sum_{\mu'} m(\mu') \times \\
&\quad \sum_{I'} \langle \tau_{\mu_1 - 1} \e_{i_1}, \ldots,
\tau_{\mu_\rho - 1} \e_{i_\rho} \mid \e^{I'}, \mu'
\rangle^{\bullet (Y_1, E)} \langle \alpha_1, \ldots, \alpha_n \mid
\e_{I'}, \mu' \rangle^{(\tilde E, E)} + R,
\end{align*}
where $R$ denotes the remaining terms which either have total
contact order smaller than $d_2$ or have number of insertions
fewer than $n$ on the $(\tilde E, E)$ side or the invariants on
$(\tilde E, E)$ are disconnected ones.

For the main terms, we claim that the total contact order $d_2 =
|\mu'|$ equals $|\mu| = \sum_{i = 1}^\rho \mu_i$. This follows
from the dimension counting on $\tilde E$ and $(\tilde E, E)$.
Indeed let $D = c_1(\tilde E).\beta + \dim \tilde E - 3$. For the
absolute invariant on $\tilde E$,
$$
\sum\nolimits_{j = 1}^n \deg \alpha_j + |\mu| - \rho +
\sum\nolimits_{j = 1}^\rho (\deg \e_{i_j} + 1) = D + n + \rho,
$$
while on $(\tilde E, E)$ (notice that now $c_1(\tilde E).\beta_2 =
d_2 c_1(\tilde E).\gamma = c_1(\tilde E).\beta$),
$$
\sum\nolimits_{j = 1}^n \deg \alpha_j + \sum\nolimits_{j =
1}^{\rho'} \deg \e_{i'_j} = D + n + \rho' - |\mu'|.
$$
Hence $(\e_I, \mu)$ appears as one of the $(\e_{I'}, \mu')$'s and
$|\mu| = |\mu'| = d_2$.

In particular, $R$ is $\T$-invariant by induction. Moreover,
$$
\deg \e_I - \deg \e_{I'} = \rho - \rho'.
$$

We will show that the highest order term in the sum consists of
the single term
$$
C(\mu)\langle \alpha_1, \ldots, \alpha_n \mid \e_I,
\mu\rangle^{(\tilde E, E)}
$$
where $C(\mu) \ne 0$.

For any $(\e_{I'}, \mu')$ in the highest order term, consider the
splitting of weighted partitions
$$
(\e_I, \mu) = \coprod\nolimits_{k = 1}^{\rho'} (\e_{I^k}, \mu^k)
$$
according to the connected components of the relative moduli of
$(Y_1, E)$, which are indexed by the contact points of $\mu'$ by
the genus zero assumption and the fact that the invariants on
$(\tilde E, E)$ are connected invariants.

Since fiber class invariants on $\mathbb{P}^1$ bundles can be
computed by pairing cohomology classes in $E$ with GW invariants
in the fiber $\mathbb{P}^1$ (c.f.\ \cite{MP}, \S1.2), we must have
$\deg \e_{I^k} + \deg \e^{i'_k} \le \dim E$ to get non-trivial
invariants. That is
$$
\deg \e_{I^k} = \sum\nolimits_{j} \deg \e_{i^k_j} \le \dim E -
\deg \e^{i'_k} \equiv \deg \e_{i'_k}
$$
for each $k$. In particular, $\deg \e_I \le \deg \e_{I'}$, hence
also $\rho \le \rho'$.

The case $\rho < \rho'$ is handled by the induction hypothesis, so
we assume that $\rho = \rho'$ and then $\deg \e_{I^k} = \deg
\e_{i'_k}$ for each $k = 1, \ldots, \rho'$. In particular $I^k \ne
\emptyset$ for each $k$. This implies that $I^k$ consists of a
single element. By reordering we may assume that $I^k = \{i_k\}$
and $(e_{I^k}, \mu^k)= \{(e_{i_k}, \mu_k)\}$.

Since the relative invariants on $Y_1$ are fiber integrals,
the virtual dimension
for each $k$ (connected component of the relative virtual moduli) is
\begin{align*}
&2\mu'_k + \dim Y_1 - 3 + 1 + (1 - \mu'_k) \\
&\qquad = (\mu_k - 1) + (\deg \e_{i_k} + 1) + (\dim E - \deg \e_{i'_k}).
\end{align*}
Together with $\deg \e_{i_k} = \deg \e_{i'_k}$ this implies that
$$
\mu'_k = \mu_k, \quad k = 1, \ldots, \rho.
$$

From the fiber class invariants consideration and
$$
\deg \e_{i_k} + \deg \e^{i'_k} = \dim E,
$$
$\e_{i_k}$ and $\e^{i'_k}$ must be Poincar\'e dual to get
non-trivial integral over $E$. That is, $\e_{i'_k} = \e_{i_k}$ for
all $k$ and $(e_{I'}, \mu') = (e_I, \mu)$. This gives the term we
expect for with $C(\mu)$ a nontrivial fiber class invariant. The
proof of Proposition~\ref{p:4.8} is complete. \medskip

The functional equations for these special absolute invariants
with descendents will be handled in \S5.

\subsection{Examples}
We consider simple $\mathbb{P}^r$ flops for $r \le 2$ in general
and for $r\ge 3$ under nefness constraint on $K_X$.

If $\beta = d\ell$, the invariant depends only on $Z$, $\alpha|_Z$
and $N_{Z/X}$. In particular
$$
\left<\alpha\right>_{g,n,d\ell}^X = \left<
p^*(\alpha|_Z)\right>_{g,n,d\ell}^{\tilde{E}}.
$$
Thus we consider $\beta \ne d\ell$. Let $\alpha_i\in H^{2l_i}(X)$.
By the divisor axiom, we may assume that $l_i \ge 2$ for all $i$.

For $\eta = (\Gamma_1, \Gamma_2, I_\rho)$ associated to
$(g,n,\beta)$, let $d$, $d_{\Gamma_1}$ and $d_{\Gamma_2}$ be the
virtual dimension (without marked points) of stable morphisms into
$X$ and relative stable morphisms into $(Y_1,E)$, $(Y_2,E)$
respectively. We have $l_1 + \cdots + l_n = d + n$. Moreover,
since $\dim E = 2r$, the degeneration formula implies that $d =
d_{\Gamma_1} + d_{\Gamma_2} - 2r\rho$.

We assume that the summand given by $\eta$ is not zero. Since
$\beta \ne d\ell$ and $A_1(Y_2)$ is spanned by $\ell$ and a fiber
line $\gamma$, we see that $\beta_1 \ne 0$ and $\Gamma_1 \ne
\emptyset$.

If $\rho = 0$ then $\Gamma_2 = \emptyset$ by connectedness, and
this gives the blow-up term
$$
\left<\tilde\alpha\right>_{g,n,\phi^*\beta}^Y.
$$
So we assume that $\rho \ne 0$. By reordering, we may assume that
in the degeneration expression $\alpha_i$ appears in the $Y_1$
part for $1\le i \le m$ and $\alpha_i$ appears in the $Y_2$ part
for $m + 1\le i \le n$. By transversality, the corresponding
relative invariant is non-trivial only if $2\le l_i \le r$ for $m
+ 1 \le i \le n$. If $r = 1$ this simply means that all
$\alpha_i$'s appear in $Y_1$. In the following we abuse the
notation by writing $|\mu|$ as $\mu$.

\begin{theorem}[Li-Ruan \cite{LiRu}]
For simple $\mathbb{P}^1$ flops of threefolds with $\beta \ne
d\ell$,
$$
\left<\alpha\right>_{g,n,\beta}^X =
\left<\tilde\alpha\right>_{g,n,\phi^*\beta}^Y = \left<\T
\alpha\right>^{X'}_{g,n,\T\beta}.
$$
That is, there are no degenerate terms and hence no analytic
continuations are needed for non-exceptional curve classes.
\end{theorem}

\Proof If $r = 1$, then $(K_X.p_*\beta_2) = 0$, $d = -(K_X.\beta)$
and
\begin{align*}
(K_Y.\beta_1) &= (\phi^* K_X.\beta_1) + (E.\beta_1) =
(K_X.\phi_*\beta_1) + \mu \\
&= (K_X.(\beta - p_*\beta_2)) + \mu = (K_X.\beta) + \mu.
\end{align*}
So
$$
d_{\Gamma_1} = -(K_Y.\beta_1) + \rho - \mu = d + (\rho - 2\mu).
$$

If $\rho \ne 0$ then $d_{\Gamma_1} < d$. Since $l_i \ge 2$, we may
assume that $\alpha_i$'s are disjoint from $Z$, hence they must
all contribute to the $Y_1$ part. This forces that $\rho = 0$ and
the result follows. \Endproof

For simple $\mathbb{P}^2$ flops, non-trivial degenerate terms do
occur even for $n \le 3$ and $g = 0$. Let $v_i := |\Gamma_i|$ be
the number of connected components.

\begin{lemma}\label{dim}
For $\tilde E = \mathbb{P}_Z(N \oplus \mathscr{O})$ of a pair
$Z\subset X$,
$$
c_1(\tilde E) = ({\rm rk}\,N + 1)E + p^*c_1(X)|_Z.
$$
\end{lemma}

\Proof Indeed, from $0 \to \mathscr{O}\to \mathscr{O}(1)\otimes
p^*(N\oplus\mathscr{O}) \to T_{\tilde E/Z} \to 0$ we get
$c_1(T_{\tilde E/Z}) = ({\rm rk}\,N + 1)E + p^*c_1(N)$, so the
formula follows from $c_1(\tilde E) = c_1(T_{\tilde E/Z}) +
p^*c_1(Z)$. \Endproof

\begin{proposition}
For simple $\mathbb{P}^2$ flops, let $n \le 3$ and $\alpha_i\in
H^{2l_i}(X)$ with $l_i \ge 2$ for $i = 1,\ldots, n$. Consider
$\beta \ne d\ell$ and an admissible triple $\eta$ with $\rho \ne
0$. Then $v_1 = \rho = \mu$, $v_2 = 1$ and $l_i = 2$ for all $i$.
\end{proposition}

\Proof Since $c_1(Y_2) = 4E$ (by Lemma \ref{dim}),
we find that
\begin{align*}
d_{\Gamma_2} &= 4(E.\beta_2) + 2v_2 + \rho - \mu \\
&= 3\mu + \rho + 2v_2.
\end{align*}
So $d_{\Gamma_2} - 4\rho = 3(\mu - \rho) + 2v_2 \ge 2$.

For one-point invariants, $l_1 = d + 1 = d_{\Gamma_1} + 3(\mu -
\rho) + 2v_2 + 1 \ge d_{\Gamma_1} + 3$. It forces that $\alpha_1$
contributes in $Y_2$, hence $l_1 = 2$ and $d = 1$. But
$d_{\Gamma_1}\ge 0$ implies that $d \ge 2$, hence a contradiction.

For two-point invariants, from $l_1 + l_2 = d + 2 = d_{\Gamma_1} +
3(\mu - \rho) + 2v_2 + 2 \ge d_{\Gamma_1} + 4$ and the fact that
$\alpha_i$ contributes to the $Y_2$ part in the degeneration
formula only if $l_i = 2$, similar argument shows that the only
non-trivial case is that $l_1 = l_2 = 2$ and both $\alpha_1$ and
$\alpha_2$ contribute in $Y_2$. Moreover the equality holds hence
that $\mu = \rho$, $v_2 = 1$ and $d_{\Gamma_1} = 0$.

We now consider three-point invariants. From
$$
l_1 + l_2 + l_3 = d + 3 = d_{\Gamma_1} + (d_{\Gamma_2} - 4\rho) +
3 \ge d_{\Gamma_1} + 5,
$$
if only $\alpha_3$ contributes to $Y_2$ then $l_1 + l_2 \ge
d_{\Gamma_1} + 3 > d_{\Gamma_1} + 2$ leads to trivial invariant.
If $\alpha_2$ and $\alpha_3$ contribute to $Y_2$, then $l_1 \ge
d_{\Gamma_1} + 1$. This leads to non-trivial invariant only if
equality holds. That is, $\mu = \rho$ and $v_2 = 1$.

The remaining case is that $l_i = 2$, $\alpha_i$ contributes in
$Y_2$ for all $i =1, 2, 3$. We have $\mu = \rho$, $v_2 = 1$, $d =
3$, $d_{\Gamma_1} = 1$, $d_{\Gamma_2} = 4\rho + 2$. \Endproof

To summarize, notice that the weighted partitions associated to
the relative invariants on the $Y_2 = \tilde E$ part are of the
form $(\mu_1, \ldots, \mu_n) = (1, \ldots, 1)$ and $\deg \alpha_i
= 2$ for all $i$, thus they are of the lowest order with fixed
$|\mu|$. They can be reduced to absolute invariants readily.

For $\beta_2 = d_1\ell + d_2\gamma$, we see that $d_2 = \mu =
\rho$ and so
$$
d_{\Gamma_2} = 4d_2 + 2
$$
is independent of $d_1$. Also $d_2$ is uniquely determined by the
cohomology insertions. The presence of degenerate terms with
degree $\beta_2$ for all large $d_1$ indicates the necessity of
analytic continuations. (c.f.\ Example \ref{ex-1}.)\medskip

The same conclusion holds for $r \ge 3$ if we impose the nefness
of $K_X$. We state the result in a slightly more general form:

\begin{proposition}
Let $\phi: Y \to X$ be the blow-up of $X$ along a smooth center
$Z$ of dimension $r$ and codimension $r' + 1$ with $K_X$ nef and
$r \le r' + 1$. Then $C_\eta \ne 0$ only if $g_1 = 0$, $v_1 = \mu
= \rho \ne 0$ and $\mu_1 \equiv 1$, $v_2 = 1$.
\end{proposition}

The proof is entirely similar and we omit it.

\section{Analytic Continuations on Local Models}  \label{s:5}

The basic strategy to calculate GW invariants on local models is
similar to \S\ref{s:3}. Here we start with one point invariants on
toric varieties. The compatibility of functional equations under
the reconstruction procedure is proved with help from operators
$\delta_H$'s which generalize $q^\ell d/dq^\ell$, the one used in
Corollary \ref{q-inv}.

\subsection{One-point functions on local models}
In this section $X$ is the local model
\[
X = \mathbb{P}_{\mathbb{P}^r}(\mathscr{O}(-1)^{\oplus(r + 1)}
\oplus \mathscr{O}).
\]
The cohomology (Chow ring) is given by
$$
H^*(X) = A^*(X) = \mathbb{Z}[h, \xi] / (h^{r + 1}, (\xi - h)^{r +
1} \xi).
$$
Since $c_1(X) = (r + 2)\xi$ is semi-positive, $X$ is a semi-Fano
toric variety.

The toric fan $\triangle(X)$ of $X$ is given by one dimensional
edges
\[
w_0, \ldots, w_{r + 1}, v_0, \ldots, v_r \in \mathbb{Z}^{r + (r +
1)}
\]
such that
\[
 w_0 + w_1 + \cdots + w_{r+1} = 0, \qquad
 v_0 + \cdots + v_r = w_0 + \cdots + w_r = -w_{r + 1}.
\]
Let $\{e_i\}_{i = 0, \ldots, r - 1}$ and $\{e'_i\}_{i = 0, \ldots,
r}$ be the basis of $\mathbb{Z}^r \times \mathbb{Z}^{r + 1}$. Then
we may pick
$$
w_i = e'_i, \quad 0 \le i \le r; \qquad w_{r + 1} = -e'_0 - \cdots
- e'_r;
$$
$$
 v_i = e_i + e'_i, \quad 0\le i \le r - 1;\qquad
 v_r = -e_0 - \cdots - e_{r - 1} + e'_r.
$$
This implies the following linear equivalence of toric divisors
$$
 D_{v_0} = D_{v_1} = \cdots = D_{v_r} =: h;\qquad
 \xi:= D_{w_{r+ 1}} =D_{w_i} + D_{v_i}, \quad i = 0, \ldots, r.
$$
Thus $D_{w_i} = \xi - h$ for all $i= 0,\ldots, r$.

\begin{remark}
In terms of the homogeneous coordinate rings, $X$ is defined by an
embedding of $(\mathbb{C}^*)^2 \hookrightarrow
(\mathbb{C}^*)^{2r+1}$, which is defined by the $2 \times (2r+1)$
matrix $M: {\rm Lie}(\mathbb{C}^*)^{2r+1} \to {\rm
Lie}(\mathbb{C}^*)^{2}$
\[
 M = \left(
 \begin{matrix}
   1 & \ldots &1 &-1 &\ldots &-1 &0 \\
   0 & \ldots &0 &1 &\ldots &1 &1
 \end{matrix}
 \right),
\]
where on the first row, there are $r$ $1$'s, $r$ $(-1)$'s. The
K\"ahler cone is spanned by $h$ and $\xi$ on $H^2(X) \cong
\mathbb{C}^2$.
\end{remark}

We start with one-point descendent invariants. The toric data
allows us to apply the known results of \cite{Givental}
\cite{LLY3} directly. Let
$$
P_\beta := \frac{\displaystyle \prod\nolimits_{\rho \in
\triangle_1(X)} \prod\nolimits_{m = -\infty}^0 (D_\rho +
mz)}{\displaystyle \prod\nolimits_{\rho \in \triangle_1(X)}
\prod\nolimits_{m = -\infty}^{(\beta.D_\rho)} (D_\rho + mz)}.
$$

\begin{lemma} \label{l:locJ}
For an effective curve class $\beta = d_1 \ell + d_2\gamma$,
\[
J_X(\beta, z^{-1}) = \frac{\displaystyle \prod_{m = - \infty}^0
(\xi - h + mz)^{r + 1}}{\displaystyle \prod_{m = 1}^{d_1} (h +
mz)^{r + 1}\prod_{m = - \infty}^{d_2 -d_1} (\xi - h + mz)^{r +
1}\prod_{m = 1}^{d_2} (\xi+ mz)}.
\]
\end{lemma}

\Proof Since $(\beta.h) = d_1$ and $(\beta.\xi) = d_2$, the right
hand side is precisely $P_\beta$. $J_X(\beta, z^{-1})$ is equal to
$\Pbeta$ without change of variables (``mirror transformation'')
due to the uniqueness theorem and the fact that $P_\beta =
O(1/z^2)$ in $1/z$ power series expansion. Indeed if $d_1 \le
d_2$,
$$
 \Pbeta = \frac{1}{(d_1!)^{r + 1}((d_2 - d_1)!)^{r + 1}
d_2!}\frac{1}{z^{d_2(r + 2)}} + \cdots,
$$
while if $d_1 > d_2$ (the key observation),
$$
 \Pbeta = (\xi - h)^{r + 1}\Big(\frac{((d_1 - d_2 - 1)!)^{r +
1}}{(d_1!)^{r + 1}d_2!}(-1)^{d_1 - d_2 -1}\frac{1}{z^{d_2(r + 2) +
r + 1}} + \cdots\Big).
$$

For more details see \cite{Givental} \cite{LLY3}. \Endproof

It also follows that a presentation of the small quantum
cohomology ring is given by Batyrev's quantum ring (cf.\
\cite{CoKa}, the proof of Proposition 11.2.17). Namely for $\q_1 =
\q^\ell$ and $\q_2 = \q^\gamma$,
$$
QH^*(X) = \mathbb{C}[h, \xi][\![\q_1, \q_2]\!] / (h^{r+1} -
\q_1(\xi - h)^{r + 1}, (\xi - h)^{r + 1}\xi - \q_2 ).
$$
Though the presentation does not provide enough information for
our purpose, it does give a first test of the invariance property.

\begin{proposition}
The map $\T_X: QH^*(X)[\q_1^{-1}] \to QH^*(X')[\q_1'^{-1}]$
defined by $\T_X h = \xi' - h'$, $\T_X \xi = \xi'$, $\T_X \q_1 =
\q_1'^{-1}$ and $\T_X \q_2 = \q_1' \q_2'$ extends to a ring
isomorphism.
\end{proposition}

\Proof Since $\T_{X'}\circ \T_X = \mbox{Id}_X$, it is enough to
check that the generators of the ideal are mapped into the
corresponding ideal in the $X'$ side:
\begin{align*}
\T_X(h^{r+1} - \q_1(\xi - h)^{r + 1}) &= (\xi' - h')^{r + 1} -
\q_1'^{-1} h'^{r + 1}\\
&= -\q_1'^{-1}(h'^{r+1} - \q_1'(\xi' - h')^{r + 1});
\end{align*}
\begin{align*}
\T_X((\xi - h)^{r + 1}\xi - \q_2) &= h'^{r + 1}\xi' - \q_1'\q_2'\\
&= (h'^{r+1} - \q_1'(\xi' - h')^{r + 1})\xi' + \q_1'((\xi' -
h')^{r + 1}\xi' - \q_2').
\end{align*}
\Endproof

Note that the virtual dimension of an $n$-point invariants in
degree $\beta = d_1\ell + d_2\gamma$ is given by $D_{n,\beta} = (r
+ 2)d_2 + 2r + n - 2$, so for a fixed set of cohomology insertions
there could be at most one $d_2$ supporting non-trivial invariants
and for the corresponding $n$-point function the summation over
$d_2$ is unnecessary.

\begin{lemma}[Quasi-linearity]\label{quasi-linear}
Let $J_X :=J_X(\q,z^{-1})$. For any $\alpha \in H^*(X)$, the one
point function $\left<\tau_k \xi\alpha\right>^X$ satisfies the
functional equation (without analytic continuation):
$$
\T\left<\tau_k \xi.\alpha\right>^X = \left<\tau_k
\T(\xi.\alpha)\right>^{X'} = \left<\tau_k
\xi'.\T\alpha\right>^{X'}.
$$
Equivalently, $\T$ is linear in $J\xi$:
$$
\T(J_X\xi.\alpha) = J_{X'}\T(\xi.\alpha) = J_{X'}\xi'.\T\alpha.
$$
\end{lemma}

\Proof The key observation on $\Pbeta$ is that if $d_2 - d_1 < 0$
then the middle factor in the denominator of $\Pbeta$ goes to the
numerator instead which has a factor $(\xi - h)^{r + 1}$. Thus it
vanishes after multiplication by $\xi$. Notice that the condition
$d_2 \ge d_1$ simply corresponds to the effectivity of $\T \beta =
-d_1\ell' + d_2(\gamma' + \ell') = (d_2 - d_1)\ell' + d_2\gamma'$.

Since $J_X = \sum_{\beta\in NE(X)} \q^\beta \Pbeta$, by the above
observation $J_X \xi.\alpha$ can be written as
$$
J_X \xi.\alpha = \sum_{d_2}\frac{1}{\displaystyle \prod_{m =
1}^{d_2} (\xi + mz)} \sum_{d_1 = 0}^{d_2} \frac{\q^{d_2\gamma}
\q^{d_1\ell}. \xi.\alpha}{\displaystyle \prod_{m = 1}^{d_1} (h +
mz)^{r + 1}\prod_{m = 1}^{d_2 - d_1} (\xi - h + mz)^{r + 1}}.
$$

Notice that since the flop is an isomorphism outside $Z =
\mathbb{P}^r \subset X$, the cohomology correspondence $\T$ is the
``identity'' one on classes $\xi.\alpha$. Namely $\T h^i = (\xi' -
h')^i$ for $i \le r$ and $\T(\xi.\alpha) = \T\xi.\T\alpha =
\xi'.\T\alpha$ for any $\alpha \in H^*(X)$, Thus
$$
\T(J_X \xi.\alpha) = \sum_{d_2}\frac{1}{\displaystyle \prod_{m =
1}^{d_2} (\xi' + mz)} \sum_{d_1 = 0}^{d_2} \frac{\q^{d_2(\gamma' +
\ell')} \q^{-d_1\ell'}.\xi'.\T\alpha}{\displaystyle \prod_{m =
1}^{d_1} (\xi' - h' + mz)^{r + 1}\prod_{m = 1}^{d_2 - d_1} (h' +
mz)^{r + 1}}.
$$
By rewriting the inner summation to be on $d_1' = d_2 - d_1 \in
\{0,\ldots, d_2\}$ we arrive at the corresponding expression of
$J_{X'}\xi'.\T\alpha$.

Since for given insertion(s) there could be at most one $d_2$
supporting non-trivial invariants, we find that $\left<\tau_k
\xi\alpha\right>$ is a finite sum and $\T\left<\tau_k
\xi.\alpha\right>  = \left<\tau_k \xi'.\T\alpha\right>$ holds
without the need of analytic continuation. \Endproof

\subsection{The functional equations in general}
Write $\beta = d_1\ell + d_2\gamma$. If $d_2 = 0$, the whole
setting on Gromov-Witten invariants goes back to quantum
corrections attached to the extremal ray $\mathbb{Z}\ell$. In
\S\ref{s:3} we had seen that while $n$-point functions with $n \ge
3$ satisfy the functional equation under $\T$ up to analytic
continuation, it is not the case for $n = 2$ or descendent
invariants with $n - 3 - k < 0$.

The results in \S\ref{s:2} and the quasi-linearity lemma are the
induction basis of our discussion on functional equations up to
analytic continuation.

For a power series $f = \sum_\beta a_\beta\, \q^\beta$ and a
divisor $H$, we define the operator
$$
\delta_H f := \sum_\beta (H.\beta) a_\beta \,\q^\beta =
\Big((H.\ell)\q^\ell\frac{\partial}{\partial \q^\ell} +
(H.\gamma)\q^\gamma\frac{\partial}{\partial \q^\gamma}\Big)f.
$$
The following lemma formalizes the argument in the proof of
Corollary \ref{q-inv}:

\begin{lemma}\label{d-lemma}
The differential operator $\delta_H$ is $\T$ equivariant. That is,
$$
\T \circ \delta_H = \delta_{\T H}\circ \T.
$$
In particular, if $\T\langle \alpha \rangle \cong \langle \T\alpha
\rangle$ then $\T \delta_H\langle \alpha \rangle \cong \delta_{\T
H}\langle \T\alpha \rangle$ too.
\end{lemma}

\Proof This follows from the fact that $\T$ preserves the
Poincar\'e pairing. In explicit terms, denote by $(x, y) =
(\q^\ell, \q^\gamma)$ and $(x', y') = (\q^{\ell'}, \q^{\gamma'})$
respectively. The transformation law $x' = x^{-1}$, $y' = xy$
leads to
$$
x\frac{\partial}{\partial x} = -x'\frac{\partial}{\partial x'} +
y'\frac{\partial}{\partial y'}; \qquad y\frac{\partial}{\partial
y} = y'\frac{\partial}{\partial y'}.
$$
Hence
\begin{align*}
\T\circ \delta_H &= (\T H. \T
\ell)\Big(-x'\frac{\partial}{\partial x'} +
y'\frac{\partial}{\partial y'}\Big) + (\T
H.\T\gamma)y'\frac{\partial}{\partial y'}\\
&= (\T H.\ell')x'\frac{\partial}{\partial x'} + (\T H.\T(\gamma +
\ell))y'\frac{\partial}{\partial y'} = \delta_{\T H}\circ\T.
\end{align*}

If $\T\langle \alpha \rangle \cong \langle \T\alpha \rangle$ then
$\T \delta_H\langle \alpha \rangle = \delta_{\T
H}\T\langle\alpha\rangle \cong \delta_{\T H}\langle \T\alpha
\rangle$. \Endproof

\begin{theorem}
Let $\langle \alpha \rangle = \langle \alpha_1, \ldots, \alpha_n
\rangle$ with $\alpha_i \in H^*(X) \cup \tau_\bullet H^*(E)$. If
$d_2 \ne 0$ then
$$
\T \langle \alpha \rangle \cong \langle \T \alpha \rangle.
$$
\end{theorem}

\Proof We will prove the theorem by induction on $d_2$ and then
$n$. This is based on the following observations: (1) By the
virtual dimension count, each set of insertions can support at
most one $d_2$. (2) Under divisor relations the degree $\beta$ is
either preserved or split into effective classes $\beta = \beta_1
+ \beta_2$, so $d_2$ is split accordingly as $d_2 = d_2^L +
d_2^R$. (3) When summing over $\beta \in NE(X)$, the splitting
terms can usually be written as the product of two generating
series with no more marked points in a manner which will be clear
in each context during the proof.

For $d_2 = 0$, since $\xi|_Z = 0$ we get trivial invariant if one
of the insertions involves $\xi$. Hence by \S\ref{s:3} the
statement in the theorem holds for $d_2 = 0$ except for the unique
case $\langle h^r, h^r\rangle$. In this case, by the divisor axiom
$$
\delta_h \langle h^r, h^r\rangle = \langle h, h^r, h^r\rangle,
$$
which satisfies the functional equation up to analytic
continuation, as we had shown before through explicit formulae
incorporated with classical defect. Thus we may base our induction
on $d_2 = 0$ with special care on this case.

Let $d_2 \ge 1$. The case $n = 1$ is contained in
Lemma~\ref{quasi-linear}, so let $n \ge 2$. We may and will make
one more assumption that $\xi$ appears in some $\alpha_i$. If not,
then there will be no descendent insertions and we may write
$$
\langle \alpha_1, \ldots, \alpha_n \rangle = \langle \alpha_1,
\ldots, \alpha_n, \xi\rangle/d_2
$$
by the divisor axiom. In the following reduction steps, each term
will either have smaller $d_2$ or with this condition being
preserved.

By reordering we may assume that $\alpha_n = \tau_s \xi a$, $s \ge
0$. Write $\alpha_1 = \tau_k h^l\xi^j$. The induction procedure is
to move divisors in $\alpha_1$ into $\alpha_n$ in the order of
$\psi$, $h$ and $\xi$. That is we use induction on the following
five numbers in the alphabetical order:
$$
(d_2, n, k, l, j).
$$

For $\psi$ we use equation $\psi_1 = -\psi_n + [D_{1|n}]^{virt}$.
If $k \ge 1$ then $j \ne 0$ and we get
\begin{align*}
\langle \tau_k h^l\xi^j, \ldots, \tau_s \xi a\rangle &= -\langle
\tau_{k - 1} h^l\xi^j, \ldots, \tau_{s + 1} \xi a\rangle \\ &\quad
+ \sum_i \langle \tau_{k - 1} h^l\xi^j, \ldots, T_i\rangle \langle
T^i, \ldots, \tau_s \xi a\rangle.
\end{align*}
For each $i$, if one of $d_2^L$ and $d_2^R$ is zero then since
both terms contain $\xi$ classes the splitting term must vanish.
So we may assume that $d_2^L < d_2$ and $d_2^R < d_2$ and these
terms are done by the induction hypothesis. By performing this
procedure to $\alpha_1, \ldots, \alpha_{n - 1}$ we may assume that
the only descendent insertion is $\alpha_n$.

For $h$, if $l \ge 2$ or $l = 1$ but $j \ne 0$ we use
\eqref{e:rec1} to get
\begin{align*}
\langle h^l\xi^j, \ldots, \tau_s \xi a\rangle &= \langle h^{l -
1}\xi^j, \ldots, \tau_s \xi ah\rangle + \delta_h\langle h^{l -
1}\xi^j, \ldots, \tau_{s + 1} \xi a\rangle \\ &\quad - \sum_i
\delta_h\langle h^{l - 1}\xi^j, \ldots, T_i\rangle \langle T^i,
\ldots, \tau_s \xi a\rangle.
\end{align*}

The only case for the splitting term to have one factor to have
the same $d_2$ and $n$ is of the form
$$
\delta_h\langle h^{l - 1}\xi^j, T_i \rangle \langle T^i, \alpha_2,
\ldots, \alpha_{n - 1}, \tau_s \xi a \rangle,
$$
where the two-point invariant has $d_2^L = 0$. But then $l - 1 <
r$ forces it to vanish. We remark here that the case $d_2^L = 0$
may still support nontrivial invariants with three or more points
if $j = 0$.

By induction (and Lemma \ref{d-lemma}) we are left with the case
$\alpha_1 = h$. The divisor axiom implies that
$$
\langle h, \ldots, \tau_s \xi a\rangle = \delta_h\langle \ldots,
\tau_s \xi a\rangle + \langle \ldots, \tau_{s - 1} \xi ah\rangle.
$$
Since both terms have one less marked points, they are done by
induction.

For $\xi$, the argument is entirely similar. For $j \ge 2$, the
divisor relation says that
\begin{align*}
\langle \xi^j, \ldots, \tau_s \xi a\rangle &= \langle\xi^{j - 1},
\ldots, \tau_s \xi^2 a \rangle + \delta_\xi \langle \xi^{j - 1},
\ldots, \tau_{s + 1} \xi a\rangle \\ &\quad - \sum_i
\delta_\xi\langle \xi^{j - 1}, \ldots, T_i\rangle \langle T^i,
\ldots, \tau_s \xi a\rangle.
\end{align*}
We then have $d_2^L < d_2$ and $d_2^R < d_2$ as before. If $j = 1$
we get
$$
\langle \xi, \ldots, \tau_s \xi a\rangle = \delta_\xi\langle
\ldots, \tau_s \xi a\rangle + \langle \ldots, \tau_{s - 1} \xi^2
a\rangle
$$
and both terms have fewer marked points. The proof is complete.
\Endproof

Practically the above inductive procedure leads to explicit
determination of GW invariants, though the computations are
somewhat tedious. For the interested readers, we list the results
for the two typical series of examples of the local model of
simple $\mathbb{P}^2$ flop.

\begin{example}\label{ex-1}
Simple $\mathbb{P}^2$ flop with $d_2 = 1$, $n = 3$. The virtual
dimension is 9. Then on $X$ ($\q_1 = \q^\ell$, $\q_2 =
\q^\gamma$),
$$
\langle h^2, h^2, h^2\xi^3 \rangle = \frac{\q_1^2}{1 + \q_1}\q_2,
\quad \langle \xi^2, \xi^2, h^2\xi^3 \rangle = (1 + \q_1)\q_2,
$$
$$
\langle h\xi, h\xi, h^2\xi^3 \rangle = \langle h\xi, \xi^2,
h^2\xi^3 \rangle = \langle h\xi, h^2, h^2\xi^3 \rangle = \langle
\xi^2, h^2, h^2\xi^3 \rangle = \q_1 \q_2.
$$
Similar formulae hold on $X'$. We compute ($\q_1' = \q^{\ell'}$,
$\q_2' = \q^{\gamma'}$)
$$
\T \langle h^2, h^2, h^2\xi^3 \rangle = \frac{\q_1'^{-2} }{1 +
\q_1^{-1}}\q_1' \q_2' = \frac{1}{1 + \q_1'}\q_2';
$$
\begin{align*}
\langle \T h^2, \T h^2, \T h^2\xi^3 \rangle &= \langle (\xi' -
h')^2, (\xi' - h')^2, \T [pt]\rangle \\
&= \langle \xi'^2, \xi'^2, [pt]\rangle + 4\langle \xi'h', \xi'h',
[pt]\rangle + \langle h'^2, h'^2, [pt]\rangle\\
&\qquad -4\langle \xi'h', h'^2, [pt]\rangle -4\langle \xi'h',
\xi'^2, [pt]\rangle + 2\langle \xi'^2, h'^2, [pt]\rangle \\
&= \Big((1 + \q_1') + 4\q_1' + \frac{\q_1'^2}{1 + \q_1'} -4\q_1'
-4\q_1' + 2\q_1'\Big)\q_2'\\
&= \Big( 1 - \q_1' + \frac{\q_1'^2}{1 + \q_1'}\Big)\q_2' =
\frac{1}{1 + \q_1'}\q_2'.
\end{align*}
Thus $\T\langle h^2, h^2, h^2\xi^3 \rangle \cong \langle \T h^2,
\T h^2, \T h^2\xi^3 \rangle$. We leave the simpler verifications
on the other five cases to the readers.
\end{example}

\begin{example}
Descendent invariants for simple $\mathbb{P}^2$ flop with $d_2 =
1$ and $n = 3$.
$$
\langle h^2, h^2, \tau_4 \xi\rangle = 3\q_1 \q_2 - 6
\frac{\q_1\q_2}{1 + \q_1}, \quad \langle \xi^2, \xi^2, \tau_4
\xi\rangle = 9\q_2 + 9\q_1 \q_2,
$$
$$
\langle h\xi, h\xi, \tau_4 \xi\rangle = \langle h^2, \xi^2, \tau_4
\xi\rangle = 3\q_2,
$$
$$
\langle h\xi, h^2, \tau_4 \xi\rangle = 0, \quad \langle h\xi,
\xi^2, \tau_4 \xi\rangle = 6\q_2 + 3\q_1 \q_2.
$$
We omit the elementary verifications on functional equations.
\end{example}

\section{Mukai flops} \label{s:6}

\subsection{Twisted Mukai flops}

Consider a flopping contraction of {\it twisted Mukai type},
namely $\psi: (X,Z) \to (\bar X, S)$ with $Z = \mathbb{P}_S(F)
\rightarrow S$, ${\rm rank}\,F = r + 1$ and $N_{Z/X} =
T^*_{Z/S}\otimes \bar\psi^* L$ for some twisting line bundle $L\in
{\rm Pic}(S)$. To construct the flop, as in the case of ordinary
flops, it is natural to consider the blow-up $\phi:Y = \mbox{Bl}_Z
X \to X$ and try to contract the exceptional set $E$ in another
fiber direction. We assume that ${\rm codim}_X Z = r \ge 2$ to
exclude trivial cases.

\begin{proposition}
Twisted Mukai flops exist.
\end{proposition}

\Proof It is well known from the case of simple Mukai flops that
the fiber $E_s$ of $E \to S$ is the degree $(1,1)$ hypersurface
$H_{1,1} \subset \mathbb{P}^r\times \mathbb{P}^r$ defined by
$$
\sum_{i = 0}^r x_i y_i = 0
$$
with $N_{E/Y}|_{E_s} = \mathscr{O}_{\mathbb{P}^r\times
\mathbb{P}^r}(-1, -1)|_{H_{1, 1}}$ (cf.\ \cite{Huyb1} or the
constructions given below). With this, the same proof as in
Proposition \ref{existence} works here. Indeed, under the same
notations, we take $C_Y$ to be a line in the $\mathbb{P}^{r - 1}$
fiber of any $E_s \cong H_{1,1} \to \mathbb{P}^r$ in the second
projection. It is clear that $(K_X.C) = 0$ and then $(K_Y.C_Y) =
-(r - 1) < 0$. To see that $C_Y$ is extremal, for $c = (H.C)$,
$$
\mathscr{O}_Y(\phi^*H + cE)|_{E_s} \cong
\mathscr{O}_{\mathbb{P}^r\times \mathbb{P}^r}(0, -c)|_{H_{1,1}}
$$
and $k\phi^*L - (\phi^*H + cE)$ is a supporting divisor for
$[C_Y]$ when $k$ is large. \Endproof

The proof suggests studying twisted Mukai flops via ordinary
flops. Indeed, we will construct the local model of it as a slice
of the ordinary flop with $F' = F^*\otimes L$. This is fundamental
throughout our later discussions.

We start with an arbitrary pair $(F, F')$ of vector bundles of
rank $r+1$ and denote the corresponding maps in the ordinary
$\mathbb{P}^r$ flop by $\Phi:\mathscr{Y}\to \mathscr{X}$,
$\Phi':\mathscr{Y}\to \mathscr{X'}$, $\Psi:\mathscr{X}\to
\bar{\mathscr{X}}$ and $\Psi':\mathscr{X}'\to \bar{\mathscr{X}}$.
Also let $g = \Psi\circ\Phi = \Psi'\circ\Phi'$. The restriction
maps on the exceptional sets are denoted by $\bar \phi$,
$\bar\phi'$, $\bar\psi$, $\bar\psi'$ and $\bar g$
respectively.\small
$$
\xymatrix{&\mathscr{E} = \mathbb{P}_S(F) \times_S
\mathbb{P}_S(F')\subset \mathscr{Y} \ar[ld]_{\Phi}
\ar[rd]^{\Phi'}\ar[dd]^g&\\
Z = \mathbb{P}_S(F)\subset \mathscr{X}\;\ar[rd]_{\Psi}& &Z' =
\mathbb{P}_S(F')\subset \mathscr{X}'\ar[ld]^{\Psi'}\\
&S\subset \bar{\mathscr{X}}}\normalsize
$$\normalsize

First suppose that there exists a non-degenerate bilinear map
$$
F\times_S F' \to \eta_S
$$
with $\eta_S \in {\rm Pic}(S)$. (This happens precisely when $F'
\cong F^* \otimes \eta_S$ for some line bundle $\eta_S$.) The map
$\mathscr{O}_{\mathbb{P}(F)}(-1) \to \bar\psi^* F$ pulls back to
$\bar\phi^*\mathscr{O}_{\mathbb{P}(F)}(-1) \to \bar g^* F$, hence
leads to a natural map
$$
\mathscr{O}_{\mathscr{E}}(-1, -1) := \bar\phi^*\mathscr{O}_{Z}(-1)
\otimes_{\mathscr{E}} \bar\phi'^*\mathscr{O}_{Z'}(-1) \to \bar
g^*(F \otimes_S F') \to \bar g^* \eta_S.
$$

Notice that the normal bundle $N_{\mathscr{E}/\mathscr{Y}}$ equals
$\mathscr{O}_{\mathscr{E}}(-1, -1)$. That is, $\mathscr{Y}$ is the
total space of $\mathscr{O}_{\mathscr{E}}(-1, -1)$. Let $p:
N_{\mathscr{E}/\mathscr{Y}} \to \mathscr{E}$ be the projection
map.

We describe two equivalent ways to construct the space $Y$. The
above linear map between line bundles induces a surjective map of
invertible sheaves which fits into an exact sequence of the form
$$
0\to N_{\mathscr{E}/\mathscr{Y}}(-E) \to
N_{\mathscr{E}/\mathscr{Y}} \to \bar g^*\eta_S \to 0
$$
for an effective divisor $E \subset \mathscr{E}$. We then take $Y
= p^{-1}(E) \subset \mathscr{Y}$ to be the collection of lines
with origins in $E$. Alternatively $Y$ is simply the irreducible
component of the inverse image of the zero section of $\bar
g^*\eta_S$ in $\mathscr{Y}$ other than the zero section
$\mathscr{E}$.

Let $X = \Phi(Y)\supset Z$, $X' = \Phi'(Y)\supset Z'$, $\bar X =
g(Y) \supset S$ with restriction maps $\phi$, $\phi'$, $\psi$,
$\psi'$. By tensoring the Euler sequence
$$
0\to\mathscr{O}_Z(-1)\to \bar\psi^*F \to \mathcal{Q}\to 0
$$
with $\mathcal{S}^* = \mathscr{O}_Z(1)$ and noticing that
$\mathcal{S}^*\otimes \mathcal{Q} \cong T_{Z/S}$, we get by
duality
$$
0\to T_{Z/S}^* \to \mathscr{O}_Z(-1)\otimes \bar\psi^*F^* \to
\mathscr{O}_Z \to 0.
$$

The inclusion maps $Z\hookrightarrow X\hookrightarrow \mathscr{X}$
leads to
$$
0\to N_{Z/X} \to N_{Z/\mathscr{X}} \to N_{X/\mathscr{X}}|_Z \to 0.
$$
Here $N_{X/\mathscr{X}}|_Z = \mathscr{O}(X)|_Z = \bar\psi^*
\mathscr{O}(\bar X)|_S$. Denote $\mathscr{O}(\bar X)|_S$ by $L$.
Recall that $N_{Z/\mathscr{X}}\cong
\mathscr{O}_{\mathbb{P}_S(F)}(-1) \otimes \bar\psi^*F'$. By
tensoring with $\bar\psi^*L^*$, we get
$$
0 \to N_{Z/X}\otimes \bar\psi^*L^* \to
\mathscr{O}_{\mathbb{P}_S(F)}(-1) \otimes \bar\psi^*(F'\otimes
L^*) \to \mathscr{O}_Z \to 0.
$$
So $F' = F^*\otimes L$ if and only if $N_{Z/X}\cong T_{Z/S}^*
\otimes \bar\psi^* L$.

This is the case for twisted Mukai flops. We have then $\eta_S
\cong L$.

\subsection{Mukai flops as limits of isomorphisms}
For Mukai flops, namely $L\cong \mathscr{O}_S$, we have $F' = F^*$
with duality pairing $F\times_S F^* \to \mathscr{O}_S$.

Consider $\pi: \mathscr{Y} \to \mathbb{C}$ via
$$
\mathscr{Y} \to \bar g^* \mathscr{O}_S = \mathscr{O}_\mathscr{E}
\cong \mathscr{E}\times \mathbb{C}
\mathop{\longrightarrow}^{\pi_2} \mathbb{C}.
$$
We get a fibration with $\mathscr{Y}_t := \pi^{-1}(t)$, which is
smooth for $t \ne 0$ and $\mathscr{Y}_0 = Y \cup \mathscr{E}$. The
intersection $E = Y\cap \mathscr{E}$ restricts to the degree
$(1,1)$ hypersurface over each fiber of $\mathscr{E} \to S$.
Indeed in each fiber the equation for $\pi$ in coordinates reads
as
$$
t = \sum_{i = 0}^r x_i y_i.
$$
From this, we also have that $\mathscr{Y}_t \cong
\mathscr{E}\backslash E$ for all $t \ne 0$ under the projection
$p$.

Let $\mathscr{X}_t$, $\mathscr{X}'_t$ and $\bar{\mathscr{X}}_t$ be
the proper transforms of $\mathscr{Y}_t$ in $\mathscr{X}$,
$\mathscr{X}'$ and $\bar{\mathscr{X}}$. For $t \ne 0$, all maps in
the diagram\small
$$
\xymatrix
{&\mathscr{Y}_t\;\ar[ld] \ar[rd]& \\
\mathscr{X}_t\;\ar[rd]& &\mathscr{X}'_t\;\ar[ld]\\
&\bar{\mathscr{X}}_t&}
$$\normalsize
are isomorphisms. For $t = 0$ this is the Mukai flop. Thus local
Mukai flops are limits of isomorphisms. More precisely, we have

\begin{theorem}
There is a projective compactification $\widehat{\mathscr{Y}} \to
\mathbb{P}^1$ which deforms the projectivized local model of Mukai
flop
$$
\widehat{\mathscr{X}}_0 = \mathbb{P}_Z(T^*_{Z/S}\oplus
\mathscr{O}) \dashrightarrow \mathbb{P}_{Z'}(T^*_{Z'/S}\oplus
\mathscr{O}) = \widehat{\mathscr{X}}'_0
$$
into isomorphisms $\widehat{\mathscr{X}}_t \cong
\widehat{\mathscr{X}}'_t \cong \mathscr{E}$ for all $t \ne 0$.

Moreover, $\widehat{\mathscr{Y}} \to \mathbb{P}^1$ is the blow-up
of $\mathscr{E}\times \mathbb{P}^1$ along $E \times \{0\}$, the
degeneration to normal cone of the pair $(\mathscr{E}, E)$ with
$E$ being the relative $(1, 1)$ divisor of
$$
\mathscr{E} = \mathbb{P}_S(F) \times_S \mathbb{P}_S(F^*)
$$
over $S$. $\widehat{\mathscr{X}}_0$, $\widehat{\mathscr{X}}'_0$
and $\widehat{\bar{\mathscr{X}}}_0$ are the contractions of
$\mathscr{E} \subset \widehat{\mathscr{Y}}_0$ along the two
rulings and the double ruling respectively.
\end{theorem}

\Proof We first consider the compactified normal bundle
$$
\bar{\mathscr{Y}} = \mathbb{P}_{\mathscr{E}} (\mathscr{O}(-1, -1)
\oplus \mathscr{O}) \dashrightarrow \mathbb{P}^1
$$
which extends the map $\pi$ by sending the infinity divisor
$\mathscr{E}_\infty \cong \mathscr{E}$ to $\infty \in
\mathbb{P}^1$.

It is clear that $E_\infty := \bar{Y} \cap \mathscr{E}_\infty$ is
the ``{\it axis}'' where $\pi$ is not defined. Indeed $E_\infty$
is the boundary divisor of every $\mathscr{Y}_t$. Thus the blow-up
$$
\xymatrix{ \widehat{\mathscr{Y}} := {\rm Bl}_{E_\infty}
\bar{\mathscr{Y}} \ar[d] \ar[rd]^{\hat{\pi}} &
\\ \bar{\mathscr{Y}}\ar@{-->}[r]^\pi & \mathbb{P}^1}
$$
resolves the indeterminacy to get a morphism $\hat{\pi}:
\widehat{\mathscr{Y}} \to \mathbb{P}^1$. $\widehat{\mathscr{Y}}_t$
is the compactification of $\mathscr{Y}_t$ by adding $E$ at
infinity, hence $\widehat{\mathscr{Y}}_t \cong \mathscr{E}$ for
all $t \ne 0$.

We then have a compactified diagram as expected:\small
$$
\xymatrix
{&\widehat{\mathscr{Y}}_t\;\ar[ld] \ar[rd]& \\
\widehat{\mathscr{X}}_t\;\ar[rd]& &\widehat{\mathscr{X}}'_t\;
\ar[ld]\\
&\widehat{\bar{\mathscr{X}}}_t&}
$$\normalsize

For $t = 0$, by the very construction of Mukai flops from the
ordinary flops, we have $\mathbb{P}_Z(T^*_{Z/X}) \cong E$. So the
compactification $\widehat{\mathscr{X}}_0$ coincides with
$\mathbb{P}_Z(T^*_{Z/S}\oplus \mathscr{O})$. Similarly
$\widehat{\mathscr{X}}'_0 \cong \mathbb{P}_{Z'}(T^*_{Z'/S}\oplus
\mathscr{O})$.

For the second statement, again by our construction
$\widehat{\mathscr{Y}}_0 = \mathscr{E} \cup \widehat{Y}$ with
$\widehat{Y}$ the total space of the $\mathbb{P}^1$ bundle
$$
\mathbb{P}_E(\mathscr{O}_E(-1, -1)\oplus \mathscr{O}) \cong
\mathbb{P}_E(\mathscr{O} \oplus \mathscr{O}_E(1, 1)).
$$
This is precisely the exceptional divisor coming from  the blow-up
$$
\widetilde{\mathscr{Y}} = {\rm Bl}_{E \times \{0\}}
\mathscr{E}\times \mathbb{P}^1.
$$
The theorem follows from an easy comparison of
$\widetilde{\mathscr{Y}}$ with $\widehat{\mathscr{Y}}$. \Endproof

In particular all interesting invariants which are continuous
under deformations are preserved. For example, the diffeomorphism
type, Hodge type and quantum cohomology rings etc.. To be more
precise, since the fiber product satisfies the base change
property and for ordinary flops the fiber product equals the graph
closure, the canonical isomorphism of Chow motives of projective
local models of Mukai flops $f:X \dashrightarrow X'$ is clearly to
be induced by the correspondence $[X\times_{\bar X} X']$, which is
the $t = 0$ fiber of the graph of $\mathscr{X} \dashrightarrow
\mathscr{X}'$:
$$
\T := [X\times_{\bar X} X'] = [\bar\Gamma_f] + [\mathscr{E}] \in
A^*(X\times X')
$$
where $[\bar\Gamma_f] \equiv Y := {\rm Bl}_Z X = {\rm Bl}_{Z'}
X'$. For global (twisted) Mukai flops we also consider the fiber
product as the proposed correspondence $\T$.

The quantum cohomologies are not just isomorphic, in fact all
quantum corrections attached to the extremal ray are zero: If not,
then the deformation invariance of Gromov-Witten invariants
implies that some extremal curve class $d\ell\in NE(X)$ survives
as an effective curve in a nearby fiber as $C_t \subset
\mathscr{X}_t \cong \mathscr{X}'_t$, then the class
$$
[C_t'] = \T_t[C_t] \sim \T d\ell = -d\ell'.
$$
is both effective and anti-effective on $\mathscr{X}'$, which is a
contradiction. (For simple Mukai flops, the invariants on $d\ell$
are zero have also been proved by Hu and Zhang \cite{HuZh} by
direct computation via localizations.)

For a global Mukai flop, the local deformation equivalence may
fail to extend to a global deformation equivalence since there are
in general obstructions to extend deformations from local to
global. (For hyper-K\"ahler manifolds or more generally Calabi-Yau
manifolds such global deformations do exists.) Nevertheless,
together with the degeneration analysis, the local deformation
equivalence do lead to global results:

\begin{theorem}
For any Mukai flop $f: X \dashrightarrow X'$ (not necessarily
being simple), $X$ is diffeomorphic to $X'$ and both have
isomorphic Chow motives, Hodge structures and full Gromov--Witten
theory (in all genera) under the correspondence $\T$. Moreover,
all quantum corrections attached to the extremal ray vanish.
\end{theorem}

\Proof The diffeomorphism is obtained by patching the local
deformation equivalence and the identity map on $X \backslash Z
\cong X' \backslash Z'$.

For Chow motives, we investigate the induced mapping on Chow
groups as in \S\ref{s:2}. For any $T$, ${\rm id}_T\times f:
T\times X \dashrightarrow T\times X'$ is also a Mukai flop, with
base $S$ being replaced by $T\times S$. Since the correspondence
$\T$ is compatible with base change, to prove that $\T^*\circ \T =
\Delta_X$, by the identity principle, we only need to show that
$\T^*\T = {\rm id}$ on $A^*(X)$ for any Mukai flop. Let $p:
X\times X' \to X$ and $p': X\times X' \to X'$ be the projections.
From
$$
\T W = p'_*(([\bar\Gamma_f] + [\mathscr{E}]).p^*W)
$$
and the property of intersection product we see that the
$\T^*\circ \T = \Delta_X$ is really a local statement which
depends only on the normal bundles $N_{Z/X}$ and $N_{Z'/X'}$. Thus
the identity follows from the case of local models. Similarly
$\T\circ \T^* = \Delta_{X'}$. So $\T$ induces an isomorphism on
Chow motives of $X$ and $X'$. The Hodge realizations leads to
equivalence of Hodge structures.

Now we treat the Gromov-Witten invariants. As in the case of
ordinary flops, we consider degeneration to normal cone $W \to
\mathbb{A}^1$ of $X$ and $W' \to \mathbb{A}^1$ of $X'$
respectively. $W_0 = Y \cup X_{\rm loc}$ with $Y = {\rm Bl}_Z X$
and $X_{\rm loc} = \mathbb{P}_Z(T^*_{Z/S}\oplus \mathscr{O})$.
Similarly $W'_0 = Y' \cup X'_{\rm loc}$ with $Y' = {\rm Bl}_{Z'}
X'$ and $X'_{\rm loc} = \mathbb{P}_{Z'}(T^*_{Z'/S}\oplus
\mathscr{O})$. By definition $Y = Y'$ and we have the induced
Mukai flop for local models $f: X_{\rm loc} \dashrightarrow
X'_{\rm loc}$.

By the degeneration formula, any Gromov-Witten invariant $\langle
\alpha \rangle^X_{g, n, \beta}$ splits into sum of products of
relative invariants of $(Y, E)$ and $(X_{\rm loc}, E)$. Now we
compare it with the similar splitting of $\langle \T\alpha
\rangle^X_{g, n, \beta}$ into $(Y, E)$ and $(X'_{\rm loc}, E)$.

Notice that most of the setting on degeneration analysis in
\S\ref{s:4} is still valid in the Mukai case. In particular, the
cohomology reduction (Proposition~\ref{coh-red}) works in the
Mukai case too.

In fact, the situation now is very simple. We match the relative
invariants on $(Y, E)$ from both sides and then we need to compare
only the cases $(X_{\rm loc}, E)$ and $(X'_{\rm loc}, E)$. But
they are deformation equivalent while the deformations leave the
boundary divisor $E$ unchanged. By the deformation invariance of
(relative) Gromov-Witten theory and the fact that $\T$ is induced
from this deformation, we find that the relative invariants are
also the same on this part. Hence we have proved
$$
\langle \alpha \rangle^X_{g, n, \beta} = \langle \T\alpha
\rangle^{X'}_{g, n, \T\beta}
$$
for any $g, n, \beta$ including descendent invariants. Namely the
full GW theory on $X$ and $X'$ are equivalent.

The statement on vanishing of GW invariants of extremal rays
follows from the previous discussion. The proof is now complete.
\Endproof

\begin{remark}
Instead of using deformation invariance of relative GW theory, we
may also proceed in the same way as the case of ordinary flops, at
least for simple Mukai flops. By Proposition \ref{reduction}, the
equivalence problem is reduced to the case of absolute invariants
and then we may use the deformation invariance of absolute GW
theory to conclude. Indeed the deformation invariance of relative
GW theory can be deduced form the absolute case and the result in
\cite{MP}.
\end{remark}

\begin{remark}
For twisted Mukai flops we take $F' = F^* \otimes L$ and $\eta_S
:= L$. The pairing $F\times_S F' \to \eta_S$ is simply $F \times_S
(F^*\otimes L) \to L$. Since
$$
\bar\phi'^*\mathscr{O}_{\mathbb{P}(F^*\otimes L)}(-1) =
\bar\phi'^*(\mathscr{O}_{\mathbb{P}(F^*)}(-1)\otimes \bar\psi'^*L)
= \bar\phi'^*\mathscr{O}_{\mathbb{P}(F^*)}(-1)\otimes \bar g^*L,
$$
the linear map $\mathscr{Y} \to \bar g^*L$ is obtained by
tensoring the corresponding map for Mukai flops with $\bar g^*L$.
The inverse image of the zero section gives $Y\cup \mathscr{E}$.
Again the proper transforms of $Y$ in various spaces give rise to
the twisted Mukai flop. The difference is that since $\bar g^*L$
is not a trivial bundle, we do not have a fibration structure
$\mathscr{Y}\to \mathbb{C}$ as before. But we still get the
equivalence of Chow motives via the fiber product.
\end{remark}

\begin{example}
To see how the extra component corrects the graph closure, we
shall carry out the detailed computations for the case of simple
Mukai flops. So $Z\cong \mathbb{P}^r$, $N_{Z/X}=T^*_Z$ and
$E\subset \mathbb{P}^r\times \mathbb{P}^r$ is the universal family
of lines in $\mathbb{P}^r$ from both sides, namely, it is the
hypersurface of bi-degree $(1,1)$. By weak Lefschetz, $H^2(E)={\rm
Pic}\,E=\mathbb{Z}x|_E\oplus \mathbb{Z}y|_E$ with $x$ and $y$ the
generators of ${\rm Pic}\,\mathbb{P}^r \times\mathbb{P}^r$ as pull
backs of $h$ and $h'$. As in the ordinary case, $N_{E/Y} =
\mathscr{O}_E(-1,-1): = \bar\phi^*\mathscr{O}_Z(-1)\otimes
\bar\phi'^*\mathscr{O}_{Z'}(-1)$.

Let $\T_0 = [\bar\Gamma_f]$. The argument to compute $\T_0$ as in
the ordinary case fails precisely when $\alpha \in A^r(X)$, so we
would like to find $\T_0[Z]$. Since $\phi^*[Z] =
j_*(c_{r-1}(\mathscr{E}))$, with $\mathscr{E}$ defined by $0\to
N_{E/Y}\to\bar\phi^*N_{Z/X}\to \mathscr{E}\to 0$, we get
\begin{align*}
&c_{r - 1}(\mathscr{E}) = \big((1 - x)^{r + 1}(1 - (x +
y))^{-1}|_E \big)_{(r - 1)}\\
&= \big((x + y)^{r-1} - C^{r + 1}_1 x(x + y)^{r - 2} + \cdots
+(-1)^{r-1}C^{r + 1}_{r - 1} x^{r - 1}\big)\big|_E\\
&= (y^{r - 1} -2y^{r - 2}x + 3y^{r - 3}x^2 + \cdots + (-1)^{r-1}r
x^{r - 1})|_E.
\end{align*}
So
$$
\T_0[Z] = \phi'_*\phi^*[Z] = (-1)^{r-1}r[Z'],
$$
which implies that $\T_0$ induces isomorphism on cohomologies over
$\mathbb{Q}$, but not over $\mathbb{Z}$.

For $0< s \le r$, since $E \sim x + y$, we have
\begin{align*}
\phi^*h^s &= j_*(c_{r-1}(\mathscr{E}).\bar
\phi^* h^s)\\
&= (y^{r-1} -2y^{r-2}x + 3y^{r-3}x^2 + \cdots + (-1)^{r-1}r
x^{r-1})|_E.x^s\\
&= (y^r - y^{r-1}x + \cdots + (-1)^{r-1}yx^{r-1} +
(-1)^r(1-r)x^r)x^s\\
&= x^s y^r - x^{s + 1}y^{r-1} + \cdots + (-1)^{r-s}x^ry^s.
\end{align*}
By symmetry this implies that $\T_0(h^s) = (-1)^{r - s}h'^s$ when
$s \ne 0$.

Let $\T = [X\times_{\bar X} X'] = \T_0 + \T_1$ with $\T_1 =
Z\times_S Z' = [\mathbb{P}^r\times \mathbb{P}^r]$. We claim that
$\T_1[Z] = (-1)^r(r + 1)[Z']$ and $\T_1h^s = 0$ for $s\ne 0$.
Indeed,
$$
\T_1[Z] = p'_*(p^{-1}[Z].[Z\times Z'])
$$
with $p$ (resp.\ $p'$) the projection of $X\times X'$ to $X$
(resp.\ $X'$). Then
$$
Z^2 = c_r(N_{Z/X}) = c_r(T^*_Z) = (-1)^r\chi(\mathbb{P}^r) =
(-1)^r(r + 1).
$$
So $\T_1[Z] = p'_*([Z\times X'].[Z\times Z']) = (-1)^r(r +
1)[Z']$. For $\T_1h^s$, notice that we may choose $W \sim Z$ with
$W\cap h^s = \emptyset$. Hence $\T_1 h^s = p'_*([h^s\times
X'].[W\times Z']) = 0$.

Thus $\T(h^s) = (-1)^{r - s}h'^s$ for $0\le s\le r$ and $\T$
induces integral isomorphisms.
\end{example}

\end{document}